\documentclass[DIV=14,letterpaper]{scrartcl}

\usepackage{tikz}
\usetikzlibrary{calc}
\usetikzlibrary{matrix}
\usepackage{amsmath,amssymb,amsfonts,amsthm,subfig} 
\usepackage{verbatim}
\usepackage{color}                    
\usepackage{hyperref}                 

\usepackage{graphicx}

\newcommand{\Redink}[1]{\color{black}{#1}}

\def \p{\partial}
\def \DD{\mathbb{D}}
\def \DDK{\mathbb{D}_K}
\newcommand{\V}[1]{\mbox{\boldmath $ #1 $}}

\newcommand{\bey}{\begin{eqnarray}}
\newcommand{\eey}{\end{eqnarray}}
\newcommand{\nn}{\nonumber}
\newcommand{\beq}{\begin{equation}}
\newcommand{\eeq}{\end{equation}}
\theoremstyle{plain}
\newtheorem{thm}{\hspace{6mm}Theorem}[section]

\newtheorem{lem}{\hspace{6mm}Lemma}[section]

\theoremstyle{definition}

\theoremstyle{remark}
\newtheorem{exam}{\hspace{6mm}Example}[section]
\newtheorem{rem}{\hspace{6mm}Remark}[section]
\newcommand{\proofend}{\mbox{ }\hfill \raisebox{.4ex}{\framebox[1ex]{}}}

\vspace{5mm}

\begin{document}

\date{}
\title{Maximum principle for the finite element solution of time dependent anisotropic diffusion problems}
\author{Xianping Li \thanks{
Department of Mathematics, the University of Central Arkansas, Conway, AR 72034,
U.S.A. ({\tt XianpingL@uca.edu})}
\and Weizhang Huang \thanks{
Department of Mathematics, the University of Kansas, Lawrence, KS 66045,
U.S.A. ({\tt huang@math.ku.edu})}
}
\maketitle

\vspace{10pt}

\begin{abstract}
Preservation of the maximum principle is studied for the combination of the linear finite element method in space and
the $\theta$-method in time for solving time dependent anisotropic diffusion problems.
It is shown that the numerical solution 
satisfies a discrete maximum principle when all element angles of the mesh measured in the metric specified
by the inverse of the diffusion matrix are nonobtuse
and the time step size is bounded below and above by bounds
proportional essentially to the square of the maximal element diameter. The lower bound requirement
can be removed when a lumped mass matrix is used. In two dimensions, the mesh and time step conditions
can be replaced by weaker Delaunay-type conditions. Numerical results 
are presented to verify the theoretical findings.
\end{abstract}



\vspace{10pt}

\section{Introduction}
\label{Sec-intro}

We are concerned with the linear finite element solution of the initial-boundary value problem (IBVP) of a linear diffusion equation,
\begin{equation}
\label{ibvp-1}
\begin{cases}
 u_t - \nabla \cdot (\mathbb{D} \, \nabla u)  =  f(\V{x},t), & \quad \mbox{ in } \quad \Omega_T = \Omega \times (0,T]
\\
u(\V{x},t) = g(\V{x},t), & \quad \mbox{ on } \quad \partial \Omega \times [0, T] \\
u(\V{x},0) = u_0(\V{x}), & \quad \mbox{ in } \quad \Omega \times \{t=0\}
\end{cases}
\end{equation}
where $\Omega \subset \mathbb{R}^d \; ({\Redink{d\ge1}})$ is a connected polygonal or polyhedral domain,
$T > 0$ is a fixed time,
$f(\V{x},t)$, $g(\V{x},t)$ and $u_0(\V{x})$ are given functions, and
$\mathbb{D}$ is the diffusion matrix.
{\Redink{We assume that $\mathbb{D}= {\Redink{\mathbb{D}(\V{x}) }}$ is a general symmetric and strictly positive
definite matrix-valued function on $\Omega_T$. It includes both isotropic and anisotropic diffusion
as special examples. In the former case, $\mathbb{D}$ takes the form $\alpha(\V{x}) I$,
where $I$ is the $d\times d$ identity matrix and $\alpha = \alpha(\V{x})$ is a scalar function.
In the latter case, on the other hand, $\mathbb{D}$ has not-all-equal eigenvalues at least
on a certain portion of $\Omega_T$.
Note that we consider only time independent $\mathbb{D}$ in this work. 
In principle, the procedure used in this work can also apply to the time dependent situation.
For that situation, however, different meshes are needed
for different time steps and the numerical solution has to be interpolated between these meshes.
Then, a conservative interpolation scheme must be employed in order for
the underlying scheme to preserve the maximum principle, non-negativity, or monotonicity.
The development of conservative interpolation schemes and their use for unstructured meshes 
is an interesting research topic in its own right (e.g., see \cite{FPPGW09}) and 
beyond the scope of the current study. To avoid this possible complexity,  we restrict our attention
to the time independent diffusion matrix in this work.
}}

Anisotropic diffusion problems arise from various areas of science and engineering including
plasma physics \cite{GL09, GLT07, GYK05, NW00, SH07, Sti92},
petroleum reservoir simulation \cite{ABBM98a,ABBM98b,CSW95,EAK01,MD06},
and image processing \cite{CS00, CSV03, KM09, MS89, PM90, Wei98}.
IBVP (\ref{ibvp-1}) is a prototype of those anisotropic diffusion problems.
It satisfies the maximum principle
{\Redink{
\beq
\max\limits_{(\V{x}, t)\in \overline{\Omega}_T} v(\V{x},t)
= \max\left \{ 0,\max\limits_{(\V{x}, t)\in \partial \Omega_T} v(\V{x},t)\right \},
\quad \forall \text{ $v$ satisfying } \quad v_t - \nabla \cdot ( \mathbb{D} \nabla v) \le 0 \text{ in }\Omega_T
\label{mp-1}
\eeq
where $\partial \Omega_T$ denotes the parabolic boundary (i.e., $\p \Omega \times \{0 < t \le T\}
\cup \Omega \times \{t = 0\}$).
}}
When a standard numerical method such as a finite element or a finite difference method
is used to solve this problem, the numerical solution may violate the maximum principe
and contain spurious oscillations.
It is of practical and theoretical importance to study when a numerical solution satisfies
a discrete maximum principle (DMP)
{\Redink{(cf. (\ref{dmp}) in Sect. \ref{Sec-cond})}}
as well as develop DMP-preserving numerical schemes.

The research topic has attracted considerable attention from researchers since 1970's and success has been
made for elliptic diffusion problems; e.g, see
\cite{BKK08, Cia70, CR73, DDS04, Hua11, KK05, KKK07, KSS09, LH10, LSS07, LS08, LuHuQi2012, 
MD06, ShYu2011, Sto82, Sto86, SF73, WaZh11, YuSh2008}.
For example, it is shown in \cite{BKK08,CR73} that for isotropic diffusion problems, the requirement of all
element angles of the mesh to be nonobtuse is sufficient for the linear finite element approximation to satisfy DMP.
In two dimensions, this nonobtuse angle condition can be replaced by a weaker, so-called Delaunay
condition \cite{SF73} which requires the sum of any pair of angles facing a common interior edge to be less than
or equal to $\pi$. For anisotropic diffusion problems, Dr\v{a}g\v{a}nescu et al. \cite{DDS04} show that
the nonobtuse angle  condition fails to guarantee DMP satisfaction for a linear finite element approximation.
Various techniques have been proposed to reduce spurious oscillations, including local matrix modification
\cite{KSS09,LS08}, mesh optimization \cite{MD06}, and mesh adaptation \cite{LSS07}.
An anisotropic nonobtuse angle condition, which uses element angles measured
in the metric specified by $\DD^{-1}$ instead of angles measured in the Euclidean metric
(as in the nonobtuse angle condition), is developed in \cite{LH10} to guarantee DMP satisfaction 
for anisotropic diffusion problems. A weaker, Delaunay-type mesh condition is obtained in \cite{Hua11}
for two-dimensional problems.
The results of \cite{Hua11,LH10} are extended in \cite{LuHuQi2012} to problems containing
convection and reaction terms.

On the other hand, less progress has been made for time-dependent problems; e.g., see
\cite{Das2006,Far2010, FH06, FH07, FH09, FHK05, FKK09, Fuj1973, Har04, LE93, MVK89, SH07, TW08, VKH08, YG06}.
Most of the existing research has focused on isotropic diffusion problems.
For example, Fujii \cite{Fuj1973} considers the heat equation and
shows that the time step size should be bounded from below and above 
for a linear finite element approximation to satisfy DMP when
the mesh satisfies the nonobtuse angle condition. He also shows that
the lower bound requirement can be removed when a lumped mass matrix is used.
{\Redink{The study is extended in \cite{Far2010} to a more general isotropic diffusion problem with a reaction term.}}
Thom{\'e}e and Wahlbin \cite{TW08} consider general anisotropic diffusion problems and show
that a semi-discrete conventional finite element solution does not satisfy DMP in general.
Slope limiters are employed in \cite{SH07} to improve DMP satisfaction for anisotropic thermal
conduction in magnetized plasmas.
Nonlinear finite volume methods are developed by Le Potier \cite{LePot05,LePot09}
for time dependent problems. 

The objective of this paper is to investigate conditions for the finite element approximation
of IBVP (\ref{ibvp-1})  to satisfy DMP for a general diffusion matrix function.
We are particularly interested in lower and upper bounds on the time step size
when the $\theta$-method and the conventional
linear finite element method are used for temporal and spatial discretization, respectively.
Two types of simplicial mesh are considered, meshes
satisfying the anisotropic nonobtuse angle condition \cite{LH10} or
a Delaunay-type mesh condition \cite{Hua11}. It is known that those meshes
lead to DMP-satisfaction linear finite element approximations
to steady-state anisotropic diffusion problems. A lumped mass matrix is also studied.
The results obtained in this paper can be viewed as a generalization of Fujii's \cite{Fuj1973} to
anisotropic diffusion problems although such generalization is not trivial. 
 
The outline of this paper is as follows. In Sect. \ref{Sec-fem}, the linear finite element solution of IBVP (\ref{ibvp-1})
is described. Sect. \ref{Sec-cond} is devoted to the development of DMP-satisfaction conditions.
Numerical examples are presented in Sect. \ref{Sec-results} to verify the theoretical findings.
Finally, Sect. \ref{Sec-con} contains conclusions.

\section{Linear finite element formulation}
\label{Sec-fem}

Consider the linear finite element solution of IBVP (\ref{ibvp-1}). 
Assume that an affine family of simplicial triangulations $\{ \mathcal{T}_h \}$ is given for the physical
domain $\Omega$. Define
\[
U_g = \{ v \in H^1(\Omega) \; | \: v|_{\p \Omega} = g\}.
\]
Denote the linear finite element space associated with mesh $\mathcal{T}_h$ by $U_g^h$.
A linear finite element solution $u^h(t) \in U_g^h$  for $t \in (0, T]$ to IBVP (\ref{ibvp-1}) is defined by
\beq
\label{fem-form}
\int_{\Omega} \frac{\p u^h}{\p t} \, v^h  d\V{x} + 
\int_{\Omega} (\nabla v^h)^T \; \mathbb{D} \nabla u^h d\V{x} =
 \int_{\Omega} f \, v^h d\V{x}, \quad \forall v^h \in U_0^h
\eeq
where $U_0^h = U_g^h$ with $g = 0$. This equation can be rewritten as
\beq
\label{disc-2}
\sum_{K \in \mathcal{T}_h} \int_{K} \frac{\p u^h}{\p t} \, v^h d\V{x} + 
\sum_{K \in \mathcal{T}_h} |K| (\nabla v^h)^{T} \, \DDK \, \nabla u^h d\V{x} =
 \sum_{K \in \mathcal{T}_h} \int_{K} f \, v^h d\V{x}, \quad \forall v^h \in U_0^h
\eeq
where $|K|$ is the volume of element $K$ and 
\[
\DDK = \frac{1}{|K|} \int_{K} \DD \, d\V{x} .
\]

Equation (\ref{disc-2}) can be expressed in a matrix form.
Denote the numbers of the elements, vertices, and interior vertices of $\mathcal{T}_h$
by $N_e$, $N_v$, and $N_{vi}$, respectively. Assume that the vertices are ordered in such a way that
the first $N_{vi}$ vertices are the interior vertices. Then $U_0^h$ and $u^h$ can be expressed as
\bey
&& U_0^h = \text{span} \{ \phi_1, \cdots, \phi_{N_{vi}} \},
\nn \\
\label{soln-approx}
&& u^h = \sum_{j=1}^{N_{vi}} u_j \phi_j + \sum_{j=N_{vi}+1}^{N_{v}} u_j \phi_j ,
\eey
where $\phi_j$ is the linear basis function associated with the $j^{\text{th}}$ vertex, $\V{a}_j$.
We approximate the boundary and {\Redink{initial conditions}} in (\ref{ibvp-1}) as
\beq
u_j(t) = g_j \equiv g(\V{a}_j, t), \quad j = N_{vi}+1, ..., N_v 
\label{fem-bc}
\eeq
{\Redink{
\beq
u_j(0) = u_0(\V{a}_j), \quad j = 1, ..., N_v .
\label{fem-ic}
\eeq
}}
Substituting (\ref{soln-approx}) into (\ref{disc-2}), taking $v^h = \phi_i$ ($i=1, ..., N_{vi}$),
and combining the resulting equations with (\ref{fem-bc}), we obtain the linear algebraic system
\beq
\label{fem-sys}
M \, \frac{d \V{u}} {d t} + A \, \V{u} = \V{f},
\eeq
where $\V{u} = (u_1,..., u_{N_{vi}}, u_{N_{vi}+1},..., u_{N_v})^T$,
$\V{f} = ( f_1, ..., f_{N_{vi}}, g_{N_{vi}+1}, ..., g_{N_v} )^T$,
\beq
\label{fem-matrix}
M = \left [\begin{array}{cc} M_{11} & M_{12} \\ 0 & 0 \end{array} \right ], \qquad
A = \left [\begin{array}{cc} A_{11} & A_{12} \\ 0 & I \end{array} \right ],
\eeq
and $I$ is the identity matrix of size $(N_{v} - N_{vi})$.
The entries of mass matrix $M$, stiffness matrix $A$, and right-hand-side vector $\V{f}$ are given by
\bey
\label{mass-matrix}
m_{ij} &=& \sum_{K \in \mathcal{T}_h} \int_K \phi_j \phi_i \, d\V{x}, 
\quad i=1, ..., N_{vi},\; j=1, ..., N_{v} \\
\label{stiffness-matrix}
a_{ij} &=& \sum_{K \in \mathcal{T}_h} |K| \;(\nabla \phi_i)^T \; \mathbb{D}_K \; \nabla \phi_j,
 \quad i=1, ..., N_{vi},\; j=1, ..., N_{v} \\
f_i &=& \sum_{K \in \mathcal{T}_h} \int_K f \phi_i \, d\V{x},
\quad i=1, ..., N_{vi}.
\eey

We use the $\theta$-method with a constant time step $\Delta t$ for time integration.
Let $\V{u}^n$ and $\V{u}^{n+1}$ be the computed solutions at the current and next time steps, respectively.
Applying the $\theta$-method to the first $N_{vi}$ equations, we get
\beq
[M_{11}\;  M_{12}] \frac{\V{u}^{n+1} - \V{u}^n}{\Delta t}  + [A_{11}\; A_{12}] 
((1-\theta) \V{u}^n + \theta \V{u}^{n+1}) = \tilde{\V{f}}^{n+\theta},
\label{time-dis-1}
\eeq
where
\[
\tilde{\V{f}}^{n+\theta} = \left [f_1(t_n + \theta \Delta t),\, ...,\, f_{N_{vi}}(t_n + \theta \Delta t)\right ]^T .
\]
For the last $N_{v}-N_{vi}$ equations (corresponding to the boundary condition), we use
\beq
u_j^{n+1} = g(\V{a}_j, t_{n+1}), \quad j = N_{vi}+1, ..., N_{v} .
\label{time-dis-2}
\eeq
Combining (\ref{time-dis-1}) and (\ref{time-dis-2}), we have
\beq
B \V{u}^{n+1} = C \V{u}^{n} + \Delta t\, \V{f}^{n+\theta},
\label{fem-2}
\eeq
where
{\Redink{
\bey
B &= & \left [\begin{array}{cc} M_{11} & M_{12} \\ 0 & I \end{array} \right ] + \theta \Delta t
\left [\begin{array}{cc} A_{11} & A_{12} \\ 0 & 0 \end{array} \right ],
\label{fem-matrix-2}\\
C & = & \left [\begin{array}{cc} M_{11} & M_{12} \\ 0 & 0 \end{array} \right ] - (1-\theta) \Delta t
\left [\begin{array}{cc} A_{11} & A_{12} \\ 0 & 0 \end{array} \right ], 
\label{fem-matrix-3}
\\
\V{f}^{n+\theta} & = & \left (  f_1(t_n + \theta \Delta t),\,
\cdots,\, f_{N_{vi}}(t_n + \theta \Delta t), \, \frac{1}{\Delta t} g(\V{a}_{N_{vi} +1}, t_{n+1})\,
\cdots,\, \frac{1}{\Delta t} g(\V{a}_{N_{v}}, t_{n+1}) \right )^T ,
\label{fem-rhs}
\\
\V{u}^0 & = & \V{u}_0 = \left (u_0(\V{a}_1), ..., u_0(\V{a}_{N_v})\right )^T.
\label{fem-ic-2}
\eey
}}
It is worth noting that the right-hand side vector, $\V{f}^{n+\theta}$, is formed from the values of the right-hand side
function $f(\V{x}, t)$ and the boundary function $g(\V{x}, t)$.
We are interested in conditions under which the scheme satisfies DMP.

\section{Conditions for DMP satisfaction}
\label{Sec-cond}

In this section we develop the conditions (on the mesh and time step size) under which
scheme (\ref{fem-2}) satisfies DMP. The main tool is a result from \cite{Sto86} which states
that the solution of a linear algebraic system satisfies DMP when the corresponding coefficient
matrix is an $M$-matrix and has nonnegative row sums. We first discuss the general dimensional case
along with the anisotropic nonobtuse angle condition developed in \cite{LH10} and then
study the two dimensional case with the Delaunay-type mesh condition developed in \cite{Hua11}.

We introduce some notation.
Consider a generic element $K \in \mathcal{T}_h$
and denote its vertices by $\V{a}_1^K,\, \V{a}_2^K,\, ...,\, \V{a}_{d+1}^K$.
Denote the face opposite to vertex $\V{a}_i^K$ (i.e., the face not having $\V{a}_i^K$ as its vertex)
by $S_i^K$ and its unit inward (pointing to $\V{a}_i^K$) normal by $\V{n}_i^K$.
The distance (or height) from vertex $\V{a}_i^K$ to face $S_i^K$ is denoted by $h_i^K$.
Define $\V{q}$-vectors as
\beq
\V{q}_i^K = \frac{\V{n}_i^K}{h_i^K} , \quad i = 1, ..., d+1.
\label{q-1}
\eeq
Obviously, we have $h_i^K = 1/{\| \V{q}_i^K\|}$.

\begin{figure}[t]
\centering
\begin{tikzpicture}[scale = 0.8]
\draw [thick] (-1,1) -- (3,0) -- (2, 2.3) -- cycle;
\draw [below] (-1,1) node {$\V{a}_3^K$};
\draw [below] (3,0) node {$\V{a}_1^K$};
\draw [above] (2,2.3) node {$\V{a}_2^K$};
\draw (1.4,0.8) node {$K$};

\draw [thick] (7,0) -- (9,2) -- (5.5, 2.5) -- cycle;
\draw [right] (9,2) node {$\widetilde{\V{a}}_3^K$};
\draw [above] (5.5,2.5) node {$\widetilde{\V{a}}_1^K$};
\draw [below] (7,0) node {$\widetilde{\V{a}}_2^K$};
\draw (7.2,1.5) node {$\widetilde{K}$};

\draw [thick] (-8,0) -- (-4,0) -- (-6, 2.3) -- cycle;
\draw [below] (-8,0) node {$\hat{\V{a}}_3$};
\draw [below] (-4,0) node {$\hat{ \V{a}}_1$};
\draw [above] (-6,2.3) node {$\hat{ \V{a}}_2$};
\draw (-6,0.8) node {$\hat K$};

\draw [thick,->] (3.5,1) -- (5.5,1);
\draw [above] (4.5,1) node {$\mathbb{D}_K^{-\frac{1}{2}} \V{x}$};

\draw [thick,->] (-3.8,1) -- (-1.8,1);
\draw [above] (-2.8,1) node {$F_K$};

\end{tikzpicture}
\caption{Sketch of coordinate transformations from $\hat K$ to $K$ and to $\widetilde{K}$.
{\Redink{Here, $\hat K$ is the reference element and $F_K$ is the affine mapping from $\hat K$ to
element $K$.}}
}
\label{f-0}
\end{figure}
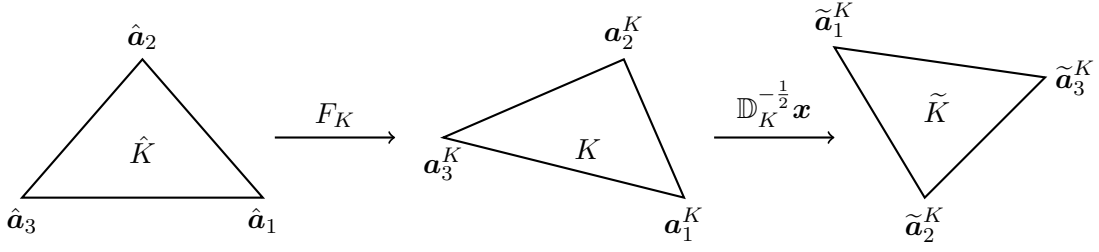

We now consider the mapping $\mathbb{D}_K^{-\frac{1}{2}} \V{x}:\, K \to \widetilde{K}$; see Fig.~\ref{f-0}.
The $\V{q}$-vectors and heights associated with $\widetilde{K}$ are denoted by
$\widetilde{\V{q}}_i^K$, and $\widetilde{h}_i^K$. We have relations
\beq
\widetilde{\V{a}}_i^K = \mathbb{D}_K^{-\frac 1 2} \V{a}_i^K, \quad
\widetilde{S}_i^K = \mathbb{D}_K^{-\frac 1 2} S_i^K, \quad
|\widetilde{K}| = \mbox{det}(\mathbb{D}_K)^{-\frac 1 2} |K|,\quad
\widetilde{\V{q}}_i^K = \mathbb{D}_K^{\frac 1 2} \V{q}_i^K , \quad
\widetilde{h}_i^K = \frac{1}{\| \widetilde{\V{q}}_i^K\|} .
\label{relation-1}
\eeq
The dihedral angle between surfaces $\widetilde{S}_i^K$ and $\widetilde{S}_j^K$ $(i \neq j)$ is denoted by
$\widetilde{\alpha}_{ij}^K$. It can be expressed as
\beq
\cos (\widetilde{\alpha}_{ij}^K) = - \frac{(\widetilde{\V{q}}_i^K)^T \widetilde{\V{q}}_j^K}
{\| \widetilde{\V{q}}_i^K \| \cdot \|\widetilde{\V{q}}_j^K \|}
= - \frac{(\V{q}_i^K)^T \mathbb{D}_K \V{q}_j^K }{\| \V{q}_i^K\|_{\mathbb{D}_K} \| \V{q}_j^K\|_{\mathbb{D}_K} },\quad
i \ne j 
\label{dihedral-2}
\eeq
where $\| \V{q}_i^K\|_{\mathbb{D}_K} = \sqrt{(\V{q}_i^K)^T \mathbb{D}_K \V{q}_i^K}$. Note that
$\widetilde{\alpha}_{ij}^K$ can be considered as {\em a dihedral angle of $K$ measured in the metric specified by
$\mathbb{D}_K^{-1}$. } 

\subsection{General dimensional case: $d \ge 1$}

We now are ready for the development of the DMP satisfaction conditions for scheme (\ref{fem-2}) for the general dimensional case.
We first have the following four lemmas.

\begin{lem}
\label{lem3.1}
For any element $K \in \mathcal{T}_h$ and $i, j = 1, ..., d+1$,
\beq
(\nabla \phi_i )^T \mathbb{D}_K \nabla \phi_j = 
\begin{cases}
- \displaystyle{\cos (\widetilde{\alpha}_{ij}^K ) \over \widetilde{h}_i^K \widetilde{h}_j^K}, & \quad \text{ for }i \neq j \\
\displaystyle{ 1 \over (\widetilde{h}_i^K)^2 }, & \quad \text{ for } i = j 
\end{cases}
\label{lem3.1-1}
\eeq
where $\phi_i$ and $\phi_j$ are the linear basis functions associated with the vertices $\V{a}_i^K$ and $\V{a}_j^K$,
respectively. In two dimensions ($d = 2$),
\beq
|K| (\nabla \phi_i )^T \DDK \nabla \phi_j = 
- \frac{\sqrt{\det(\DDK)}}{2} \cot (\widetilde{\alpha}_{ij}^K ), \quad i \ne j, \quad i, j = 1, 2, 3. 
\label{lem3.1-2}
\eeq
\end{lem} 

{\bf Proof.} see \cite{Hua11,LuHuQi2012}.
\proofend

\begin{lem}
\label{lem3.2}
The stiffness matrix $A$ defined in (\ref{fem-matrix}) and (\ref{stiffness-matrix})
is an $M$-matrix and has nonnegative row
sums if the mesh satisfies the anisotropic nonobtuse angle condition
\beq
0 < \widetilde{\alpha}_{ij}^K \le \frac{\pi}{2},\quad \forall i, j = 1, ..., d+1, i\ne j, \; \forall K \in \mathcal{T}_h .
\label{anoac-1}  
\eeq
\end{lem}

{\bf Proof.} See \cite[Theorem 2.1 and its proof]{LH10}.
\proofend

\begin{lem}
\label{lem3.3}
Matrix $B$ defined in (\ref{fem-matrix-2}) ($0 < \theta \le 1$) is an $M$-matrix if the mesh satisfies (\ref{anoac-1})
and the time step size satisfies
\beq
\Delta t \ge \frac{1}{\theta (d+1)(d+2)}
\max\limits_{K \in \mathcal{T}_h}
\max\limits_{\substack{i, j = 1, ..., d+1\\ i \neq j}}
\frac{h_i^K h_j^K}{\cos(\widetilde{\alpha}_{ij}^K)\lambda_{min}(\mathbb{D}_K)} .
\label{dt-1}
\eeq
\end{lem}

{\bf Proof.}
We first show that $M + \theta \Delta t A$ is a $Z$-matrix, i.e., it has positive diagonal and nonpositive
off-diagonal entries. From (\ref{fem-matrix}) we only need to show
\bey
&& m_{i i}+\theta \Delta t \; a_{i i} > 0, \qquad i = 1, ..., N_{vi}
\label{cond-$M$-matrix-1}
\\
&& m_{ij}+\theta \Delta t \; a_{ij} \le 0 \qquad \forall i \ne j,\; i = 1, ..., N_{vi}, \; j = 1, ..., N_{v}  .
\label{cond-$M$-matrix-2}
\eey
Let $\omega_i$ be the patch of the elements containing vertex $\V{a}_i$. Notice that $\nabla \phi_i = 0$
when $K \notin \omega_i$. Recall from \cite{Cia78} that 
\beq
\label{int_phi}
\int_{K} \phi_i \phi_j \, d\V{x} = \frac{|K|}{(d+1)(d+2)}, \qquad \int_{K} \phi_i^2 \, d\V{x} = \frac{2|K|}{(d+1)(d+2)}.
\eeq
Then (\ref{cond-$M$-matrix-1}) follows immediately from (\ref{mass-matrix}) and Lemma~\ref{lem3.2}.

For  (\ref{cond-$M$-matrix-2}), from (\ref{mass-matrix}), (\ref{stiffness-matrix}), and (\ref{int_phi}) we have
\bey
\nn
m_{ij} + \theta \Delta t\, a_{ij} & = & 
\sum_{K \in \mathcal{T}_h} \int_K \phi_j \phi_i \, d\V{x}
+ \theta \Delta t \, \sum_{K \in \mathcal{T}_h} |K| \;(\nabla \phi_i)^T \; \mathbb{D}_K \; \nabla \phi_j \\
\nn
& = & \sum_{K \in \omega_i \cap \omega_j} \left ( \int_K \phi_j \phi_i \, d\V{x}
+ \theta \Delta t \, |K| \;(\nabla \phi_i)^T \; \mathbb{D}_K \; \nabla \phi_j \right)
\nn \\
& = & \sum_{K \in \omega_i \cap \omega_j} \left ( \int_K \phi_j \phi_i \, d\V{x}
+ \theta \Delta t \, |K| \;(\nabla \phi_{i_K})^T \; \mathbb{D}_K \; \nabla \phi_{j_K} \right) ,
\label{B-ij}
\eey
where $i_K$ and $j_K$ denote the local indices (on element $K$) of vertices $\V{a}_i$ and $\V{a}_j$.
From (\ref{int_phi}) and Lemma~\ref{lem3.1}, we get
\beq
m_{ij} + \theta \Delta t\, a_{ij}
= \sum_{K \in \omega_i\cap \omega_j} |K| \left ( \frac{1}{(d+1)(d+2)} 
- \theta \Delta t \frac{\cos(\widetilde{\alpha}_{i_K j_K}^K)}{\widetilde{h}_{i_K}^K \widetilde{h}_{j_K}^K}  \right )  .
\eeq
The right-hand side term is nonpositive if
\beq
\label{cond-dt-1}
\Delta t \ge \frac{1}{\theta (d+1)(d+2)}
\max\limits_{K \in \mathcal{T}_h}
\max\limits_{\substack{i, j = 1, ..., d+1\\ i \neq j}}
\frac{\widetilde{h}_i^K \widetilde{h}_j^K}{\cos(\widetilde{\alpha}_{ij}^K)} .
\eeq
Moreover, (\ref{relation-1}) implies
\[
\widetilde{h}_{i}^K = \frac{1}{\| \widetilde{\V{q}}_i \|} = \frac{1}{\sqrt{\V{q}_i^T \mathbb{D}_K \V{q}_i }} .
\]
Thus, we have
\beq
\frac{h_i^K}{\sqrt{\lambda_{max}(\mathbb{D}_K)}} \le
\widetilde{h}_{i}^K \le \frac{h_i^K}{\sqrt{\lambda_{min}(\mathbb{D}_K)}} .
\label{h-1}
\eeq
From this, we can see that (\ref{dt-1}) implies (\ref{cond-dt-1}). Hence, we have shown
that $B$ is a $Z$-matrix when (\ref{dt-1}) holds.

To show $B$ is an $M$-matrix, we recall from (\ref{fem-matrix-2}) that
\[
B = \left [\begin{array}{cc} M_{11} + \theta \Delta t  A_{11} & M_{12} + \theta \Delta t A_{12}
\\ 0 & I \end{array} \right ] .
\]
The fact that $B$ is a $Z$-matrix means that $M_{11} + \theta \Delta t  A_{11}$ is also a $Z$-matrix
and $M_{12} + \theta \Delta t A_{12} \le 0$.
It is easy to show that $M_{11} + \theta \Delta t  A_{11}$ is positive definite, which in turn implies
$M_{11} + \theta \Delta t  A_{11}$ is an $M$-matrix. Notice
\[
B^{-1} = \left [\begin{array}{cc} (M_{11} + \theta \Delta t  A_{11})^{-1} &
-  (M_{11} + \theta \Delta t  A_{11})^{-1} (M_{12} + \theta \Delta t A_{12})
\\ 0 &  I \end{array} \right ] .
\]
This means $B^{-1} \ge 0$ and hence $B$ is an $M$-matrix.
\proofend

\begin{lem}
\label{lem3.4}
Matrix $C$ defined in (\ref{fem-matrix-3}) ($0 \le \theta \le 1$) is nonnegative if the mesh satisfies (\ref{anoac-1})
and the time step size satisfies
\beq
\Delta t \le \frac{2}{(1-\theta) (d+1)(d+2) } \min\limits_{K \in \mathcal{T}_h}
\min\limits_{i= 1, ..., d+1} \frac{(h_i^K)^2}{\lambda_{max}(\mathbb{D}_K)}.
\label{dt-2}
\eeq
\end{lem}

{\bf Proof.}
For off-diagonal entries ($i\ne j,\; i=1, ..., N_{vi},\; j = 1, ..., N_{v}$), $m_{ij} - (1-\theta) \Delta t \, a_{ij}$, are
nonnegative since $a_{ij}\le 0$ under condition (\ref{anoac-1}) (cf. Lemma~\ref{lem3.2}) and $m_{ij} \ge 0$
from definition (\ref{mass-matrix}). To see if the diagonal entries are also nonnegative, from 
(\ref{mass-matrix}), (\ref{stiffness-matrix}),  and (\ref{int_phi}) we have
\bey
\label{cond-Positive-3}
m_{ii} - (1-\theta) \Delta t \, a_{ii} 
& = & \sum_{K \in \omega_i} |K| \left ( \frac{2}{(d+1)(d+2)} 
- \frac{(1-\theta) \Delta t }{(\widetilde{h}_{i_K}^K)^2}  \right ). 
\eey
The right-hand side term is nonnegative if
\[
\Delta t \le \frac{2}{(1-\theta) (d+1)(d+2) } \min\limits_{K \in \mathcal{T}_h}
\min\limits_{i= 1, ..., d+1}
(\widetilde{h}_{i}^K)^2 .
\]
From (\ref{h-1}) we see that this condition holds when (\ref{dt-2}) is satisfied.
\proofend

\vspace{20pt}

We are now in a position to prove our first main theoretical result.

\begin{thm}
\label{thm3.1}
Scheme (\ref{fem-2}) satisfies a discrete maximum principle if the mesh satisfies
the anisotropic nonobtuse angle condition (\ref{anoac-1}) and the time step size satisfies
(\ref{dt-1}) and (\ref{dt-2}), i.e.,
\bey
&& \frac{1}{\theta (d+1)(d+2)}
\max\limits_{K \in \mathcal{T}_h}
\max\limits_{\substack{i, j = 1, ..., d+1\\ i \neq j}}
\frac{h_i^K h_j^K}{\cos(\widetilde{\alpha}_{ij}^K)\lambda_{min}(\mathbb{D}_K)}
\nn \\
&& \qquad \qquad \qquad
\quad \le \; \Delta t \; \le \frac{2}{(1-\theta) (d+1)(d+2) } \min\limits_{K \in \mathcal{T}_h}
\min\limits_{i= 1, ..., d+1} \frac{(h_i^K)^2}{\lambda_{max}(\mathbb{D}_K)} .
\label{dt-3}
\eey
\end{thm}

{\bf Proof.}
Scheme (\ref{fem-2}) can be expressed as
\beq
\mathbb{A} \mathbb{U} = \mathbb{F},
\label{fem-3}
\eeq
where
\beq
\mathbb{A} = \left [
\begin{array}{cccc}
I & 0 & & \\
- C & B   &    &  \\
& - C & B  & \\
 &  & \cdots & \cdots \\
& &  - C & B
\end{array}
\right ] ,
\qquad
\mathbb{U} = \left [ \begin{array}{c} \V{u}^0 \\
\V{u}^1 \\ \V{u}^2 \\ \vdots  \\ \V{u}^N \end{array} \right ] ,
\qquad \mathbb{F} =   \left [ \begin{array}{c}  \V{u}_0 \\ \Delta t \V{f}^{\theta}
\\ \Delta t \V{f}^{1+\theta} \\ \vdots \\  \Delta t \V{f}^{N-1+\theta} \end{array} \right ] ,
\label{big-A}
\eeq
and $B$ and $C$ are defined in (\ref{fem-matrix-3}).
Scheme (\ref{fem-3}) satisfies a DMP if
coefficient matrix $\mathbb{A}$ is an $M$-matrix and has nonnegative row sums.
From Lemmas \ref{lem3.3} and \ref{lem3.4} we know that $B$
is an $M$-matrix and $C \ge 0$. As a result, $\mathbb{A}$ is a $Z$-matrix.
Moreover, we can show $\mathbb{A}^{-1} \ge 0$. Indeed, from (\ref{fem-3}) we know that
{\Redink{$\V{u}^0 = \V{u}_0$ and thus if $\V{u}_0 \ge 0$, we have $\V{u}^0 \ge 0$.}}
Next, from the scheme we have $\V{u}^1 = B^{-1} \Delta t \V{f}^{\theta}
+ B^{-1} C \V{u}^0$. Recall that $C\ge 0$ and $B$ is an $M$-matrix and thus $B^{-1}\ge 0$.
Combining these results, we can conclude that
$\V{f}^{\theta} \ge 0$ implies $\V{u}^1 \ge 0$. Similarly, 
we can show $\V{u}^n\ge 0$ if $\V{f}^{n-1+\theta} \ge 0$, $n = 2, ..., N$.
Thus, we have shown that $\mathbb{F} \ge 0$ implies $\mathbb{U} \ge 0$.
This implies $\mathbb{A}^{-1} \ge 0$
and $\mathbb{A}$ is an $M$-matrix.

We notice that the sum of each of the second to the last (block) rows is 
\[
B - C = \left [\begin{array}{cc} \Delta t \, A_{11} & \Delta t\, A_{12} \\ 0 & I \end{array} \right ] .
\]
Since $A$ has nonnegative row sums (cf. Lemma~\ref{lem3.2}), $\mathbb{A}$ has nonnegative row sums.
{\Redink{
Thus, we have proven that $\mathbb{A}$ is an $M$-matrix and has nonnegative row sums.

Form \cite{Sto86}[Theorem 1], we conclude that the solution of  (\ref{fem-3}) satisfies
\beq
\max\limits_{i=1, ..., (N+1) N_{v}} \mathbb{U}_i = \max\left \{0,
\max\limits_{i \in S(\mathbb{F}^+)} \mathbb{U}_i \right \} ,
\label{sto-1}
\eeq
where $S(\mathbb{F}^+)$ is the set of the indices with $\mathbb{F}_i > 0$. 
When $f(\V{x},t) \le 0$, from (\ref{fem-rhs}) we know that $\mathbb{F}_i > 0$ holds
only for those indices corresponding to the boundary points on $\p \Omega_T$.
Moreover, from (\ref{fem-matrix-2}), (\ref{fem-matrix-3}), and (\ref{big-A}) we see that
at the boundary points, $\mathbb{U}_i$ is equal to either the boundary function $g$ or
the initial function $u_0$. Since a piecewise linear function attains its maximum value at vertices,
(\ref{sto-1}) implies that when $f(\V{x}, t) \le 0$, the solution of (\ref{fem-2}) satisfies a DMP
\beq
\max\limits_{n = 0, ..., N} \max\limits_{\V{x} \in \overline{\Omega}} U^n(\V{x}) 
= \max\left \{ 0, \;  \max\limits_{n = 1, ..., N} \max\limits_{\V{x} \in \p \Omega} U^n(\V{x}),
\; \max\limits_{\V{x} \in \overline{\Omega}} U^0(\V{x})  \right \} ,
\label{dmp}
\eeq
where
\[
U^n(\V{x}) = \sum_{j=1}^{N_{vi}} u_j^n \phi_j(\V{x}) + \sum_{j=N_{vi}+1}^{N_{v}} u_j^n \phi_j(\V{x}),
\quad n = 0, ...., N.
\]
Hence, we have proven that scheme (\ref{fem-2}) satisfies DMP.
}}
\proofend

\vspace{10pt}

{\Redink{
\begin{rem}
\label{rem3.1}
Consider a special case with $\mathbb{D} = \alpha I$, where $\alpha$ is a positive constant.
It is known (e.g., see Emert and Nelson \cite{EN97}) that the height (or altitude),
volume, and cosine of the dihedral angles of a regular $d$-dimensional simplex $K$ are given by
\beq
h_K = e_K \sqrt{\frac{d+1}{2 d}}, 
\quad |K| = \frac{\sqrt{d+1}}{d! \sqrt{2^d}} e_K^d , \quad
\cos(\alpha_{ij}^K) = \frac{1}{d} ,
\label{simplex-1}
\eeq
where $e_K$ is the edge length.
Thus, if the elements of $\mathcal{T}_h$ are all regular simplexes, (\ref{dt-3}) reduces to
\beq
\frac{\max\limits_{K \in \mathcal{T}_h} e_K^2 }{2 \theta \alpha (d+2)}
 \le \Delta t
\le \frac{\min\limits_{K \in \mathcal{T}_h} e_K^2}{(1-\theta) \alpha d (d+2) } .
\label{dt-3+1}
\eeq
If further the mesh is uniform (and thus all mesh elements have the same volume and same edge length ($e$)),
the above condition becomes
\beq
 \frac{ e^2}{2 \theta \alpha (d+2)}
\le \Delta t  \le \frac{e^2}{(1-\theta) \alpha d (d+2) } ,
\label{dt-3+2}
\eeq
which is exactly the result of Theorem 20 of \cite{Far2010} where the maximum principle of linear finite element
approximation of isotropic diffusion problems is studied. Interestingly, we can rewrite (\ref{dt-3+2}) in terms of
the number of the elements, $N_e$. Indeed, since the mesh is uniform, the elements have a constant volume
$|\Omega|/N_e$. From (\ref{simplex-1}), we have
\[
e = \sqrt{2} N_e^{- \frac{1}{d}} \left (\frac{|\Omega| d!  }{\sqrt{d+1}}\right )^{\frac{1}{d}} .
\]
Inserting this into (\ref{dt-3+2}), we get
\beq
 \frac{ N_e^{- \frac{2}{d}}}{ \theta \alpha (d+2)} \left (\frac{|\Omega| d!  }{\sqrt{d+1}}\right )^{\frac{2}{d}}
\le \Delta t  \le \frac{2 N_e^{- \frac{2}{d}} }{(1-\theta) \alpha d (d+2) } \left (\frac{|\Omega| d!  }{\sqrt{d+1}}\right )^{\frac{2}{d}} .
\label{dt-3+3}
\eeq
\proofend
\end{rem}
}}

{\Redink{
\begin{rem}
\label{rem3.2}
Another special case is that the mesh is uniform in the metric specified by $\mathbb{D}^{-1}$.
It is known \cite{HR11} that such a mesh satisfies the so-called alignment and equidistribution conditions
\bey
&& \frac{\frac{1}{d} \text{tr} ( (F_K')^T\mathbb{D}_K^{-1}F_K') }{\text{det} ( (F_K')^T\mathbb{D}_K^{-1}F_K')^{\frac{1}{d}}}
= 1, \quad \forall K \in \mathcal{T}_h
\label{ali-1}
\\
&& |K| \sqrt{\text{det}(\mathbb{D}_K^{-1})} = \frac{\sigma_h}{N_e},
\quad \forall K \in \mathcal{T}_h
\label{eq-1}
\eey
where $\text{tr}(\cdot)$ and $\text{det}(\cdot)$ denote the trace and determinant of a matrix, $F_K'$
is the Jacobian matrix of the affine mapping $F_K$ from the reference element $\hat K$ to element $K$, and
\beq
\sigma_h = \sum\limits_{K \in \mathcal{T}_h} |K| \sqrt{\text{det}(\mathbb{D}_K^{-1})}.
\label{sigma-1}
\eeq
Geometrically, the alignment condition (\ref{ali-1}) implies that the element $\widetilde{K}$ in Fig.~\ref{f-0}
is a regular simplex while the equidistribution condition indicates that all elements have a constant
volume $\sigma_h/N_e$ in the metric $\mathbb{D}^{-1}$. 

For such a mesh, it is more suitable to replace (\ref{dt-3}) by 
\bey
&& \frac{1}{\theta (d+1)(d+2)}
\max\limits_{K \in \mathcal{T}_h}
\max\limits_{\substack{i, j = 1, ..., d+1\\ i \neq j}}
\frac{\widetilde{h}_i^K \widetilde{h}_j^K}{\cos(\widetilde{\alpha}_{ij}^K)}
\nn \\
&& \qquad \qquad \qquad
\quad \le \; \Delta t \; \le \frac{2}{(1-\theta) (d+1)(d+2) } \min\limits_{K \in \mathcal{T}_h}
\min\limits_{i= 1, ..., d+1} (\widetilde{h}_i^K)^2 .
\label{dt-3+4}
\eey
Using the same procedure as in Remark~\ref{rem3.1} and noticing that $\widetilde{K}$ is regular, we can get
\beq
 \frac{ N_e^{- \frac{2}{d}}}{ \theta (d+2)} \left (\frac{\sigma_h d!  }{\sqrt{d+1}}\right )^{\frac{2}{d}}
\le \Delta t  \le \frac{2 N_e^{- \frac{2}{d}} }{(1-\theta) d (d+2) } \left (\frac{\sigma_h d!  }{\sqrt{d+1}}\right )^{\frac{2}{d}} .
\label{dt-3+5}
\eeq
Notice that the difference between (\ref{dt-3+3}) and (\ref{dt-3+5}) lies in that the factor, $|\Omega|/\alpha$,
has been replaced by the volume of $\Omega$ in the metric $\mathbb{D}^{-1}$, $\sigma_h$.
\proofend
\end{rem}
}}

\begin{rem}
\label{rem3.3}
It is known \cite{LH10} that a mesh,  generated as a uniform mesh
in the metric specified by $M_K = \theta_K \mathbb{D}_K^{-1}$ for all $K \in \mathcal{T}_h$, where $\theta_K$
is an arbitrary piecewise constant, scalar function defined on $\Omega$, satisfies the anisotropic
nonobtuse angle condition (\ref{anoac-1}). The reader is referred to \cite{LH10} for more information on
the generation of such meshes.
\proofend
\end{rem}

\vspace{10pt}

The lower bound requirement on $\Delta t$ in (\ref{dt-3}) can be avoided by using a lumped mass matrix.
In this case, scheme (\ref{fem-2}) is modified into
\begin{align}
& \left (\left [\begin{array}{cc} \bar{M}_{11} & 0 \\ 0 & I \end{array} \right ] + \theta \Delta t
\left [\begin{array}{cc} A_{11} & A_{12} \\ 0 & 0 \end{array} \right ]\right )
\V{u}^{n+1} = 
\notag \\
& \qquad \left ( \left [\begin{array}{cc} \bar{M}_{11} & 0 \\ 0 & 0 \end{array} \right ] - (1-\theta) \Delta t
\left [\begin{array}{cc} A_{11} & A_{12} \\ 0 & 0 \end{array} \right ]\right ) \V{u}^n
+ \Delta t \V{f}^{n+\theta} ,
\label{fem-4}
\end{align}
where $\bar{M}_{11}$ is the lumped mass matrix with diagonal entries 
\[
\bar{m}_{ii} =  \sum\limits_{j=1}^{N_v} m_{ij}, \quad i=1,..., N_{vi} .
\]
The following theorem can be proven in a similar manner as for Theorem \ref{thm3.1}.

\begin{thm}
\label{thm3.2}
Scheme (\ref{fem-4}) with a lumped mass matrix
satisfies a discrete maximum principle if the mesh satisfies
the anisotropic nonobtuse angle condition (\ref{anoac-1}) and the time step size satisfies
\bey
\Delta t \le \frac{1}{(1-\theta) (d+1) } \min\limits_{K \in \mathcal{T}_h}
\min\limits_{i= 1, ..., d+1} \frac{(h_i^K)^2}{\lambda_{max}(\mathbb{D}_K)} .
\label{dt-4}
\eey
\end{thm}

{\Redink{
\begin{rem}
\label{rem3.4}
If the mesh is uniform in the metric specified by $\mathbb{D}^{-1}$, the condition (\ref{dt-4}) reduces to
\beq
 \Delta t  \le \frac{N_e^{- \frac{2}{d}} }{(1-\theta) d } \left (\frac{\sigma_h d!  }{\sqrt{d+1}}\right )^{\frac{2}{d}} ,
\label{dt-4+1}
\eeq
where $\sigma_h$ is defined in (\ref{sigma-1}).
\proofend
\end{rem}
}}

\subsection{Two dimensional case: $d=2$}

{\Redink{The results in the previous subsection are valid for all dimensions. However, }}
it is known \cite{Hua11} that
a Delaunay-type mesh condition, which is weaker than the nonobtuse angle condition (\ref{anoac-1}),
is sufficient for a linear finite element approximation to satisfy DMP in two dimensions 
for steady-state problems. It is interesting to know if this is also true for time-dependent problems.

Consider an arbitrary interior edge $e_{ij}$. Denote the two vertices of the edge by $\V{a}_i$ and $\V{a}_j$
and the two elements sharing this common edge by $K$ and $K'$. Let the local indices of the vertices on $K$ be
$i_K$ and $j_K$. The angle of $K$ opposite $e_{ij}$ is denoted by $\alpha_{{i_K}, {j_K}}^K$ (when measured in
the Euclidean metric) and by $\widetilde{\alpha}_{{i_K}, {j_K}}^K$ when measured in the metric $\DDK^{-1}$.
Similarly, we have $\alpha_{{i_{K'}}, {j_{K'}}}^{K'}$ and $\widetilde{\alpha}_{{i_{K'}}, {j_{K'}}}^{K'}$.

\begin{lem}
\label{lem3.5}
The stiffness matrix $A$ defined in (\ref{fem-matrix}) and (\ref{stiffness-matrix}) is an $M$-matrix and has nonnegative row sums if the mesh satisfies the Delaunay-type mesh condition
\bey
&& \frac{1}{2} \left[ \frac{}{}
\widetilde{\alpha}_{{i_K}, {j_K}}^K 
+ arccot \left( \sqrt{\frac{\det(\DDK)}{\det(\mathbb{D}_{K'})}} \cot(\widetilde{\alpha}_{{i_K}, {j_K}}^K) \right)
\right .
\nn \\
&& \qquad \left. \frac{}{} 
+ \widetilde{\alpha}_{{i_{K'}}, {j_{K'}}}^{K'} 
+ arccot \left( \sqrt{\frac{\det(\mathbb{D}_{K'})}{\det(\DDK)}} \cot(\widetilde{\alpha}_{{i_{K'}}, {j_{K'}}}^{K'} ) \right)
\right]  \le \pi, \qquad \forall\; \text{interior edges } e_{ij} .
\label{delaunay-1}  
\eey
\end{lem}

{\bf Proof.} See \cite[Theorem 4.1]{Hua11}.
\proofend

\begin{lem}
\label{lem3.6}
Matrix $B$ defined in (\ref{fem-matrix-2}) ($0 < \theta \le 1$) is an $M$-matrix if the mesh satisfies (\ref{delaunay-1})
and the time step size satisfies
\beq
\Delta t \ge \frac{1}{6 \theta}
\max\limits_{e_{ij}}
\frac{|K|+|K'|}{\sqrt{\det(\DDK)}\cot(\widetilde{\alpha}_{{i_K}, {j_K}}^K)
+\sqrt{\det(\mathbb{D}_{K'})}\cot(\widetilde{\alpha}_{{i_{K'}}, {j_{K'}}}^{K'})} ,
\label{dt2-1}
\eeq
where the maximum is taken over all interior edges and
$K$ and $K'$ are the two elements sharing the common edge $e_{ij}$.
\end{lem}

{\bf Proof.} Inequality (\ref{dt2-1}) follows from (\ref{lem3.1-2}), (\ref{int_phi}), and (\ref{B-ij}). 
\proofend

\begin{lem}
\label{lem3.7}
Matrix $C$ defined in (\ref{fem-matrix-3}) ($0 < \theta \le 1$) is nonnegative if the mesh satisfies (\ref{delaunay-1})
and the time step size satisfies
\beq
\Delta t \le \frac{1}{6 (1-\theta)} \min\limits_{i}
\frac{|\omega_i|}{\sum\limits_{K\in \omega_i} |K|\, \lambda_{max}(\DDK)\, (h_{i_K}^K)^{-2} } ,
\label{lem3.7-1}
\eeq
where the minimum is taken over all interior vertices and $\omega_i$ is the patch of the elements containing
$\V{a}_i$ as its vertex.
\end{lem}

{\bf Proof.} The proof is similar to that of Lemma~\ref{lem3.4}. Indeed, Lemma~\ref{lem3.5} implies that
the off-diagonal entries of $C$ are nonnegative under condition (\ref{delaunay-1}). For diagonal entries,
from (\ref{cond-Positive-3}) we get
\[
m_{ii} - (1-\theta)\Delta t a_{ii} = 
\frac{ |\omega_i|}{6} - (1-\theta) \Delta t \sum\limits_{K \in \omega_i} 
\frac{|K|}{(\widetilde{h}_{i_K}^K)^2} .
\]
From (\ref{h-1}), we can see that the right-side term of the above equation is nonnegative when (\ref{lem3.7-1}) holds.
\proofend

\vspace{10pt}

Using the above results we can prove the following theorems in a similar manner
as for Theorems~\ref{thm3.1} and \ref{thm3.2}.

\begin{thm}
\label{thm3.3}
In two dimensions, scheme (\ref{fem-2}) satisfies a discrete maximum principle if the mesh satisfies
the Delaunay-type mesh condition (\ref{delaunay-1}) and the time step size satisfies
(\ref{dt2-1}) and (\ref{lem3.7-1}), i.e.,
\bey
&& \frac{1}{6 \theta}
\max\limits_{e_{ij}}
\frac{|K|+|K'|}{\sqrt{\det(\DDK)}\cot(\widetilde{\alpha}_{{i_K}, {j_K}}^K)
+\sqrt{\det(\mathbb{D}_{K'})}\cot(\widetilde{\alpha}_{{i_{K'}}, {j_{K'}}}^{K'})}
\nn \\
&& \qquad \qquad \qquad 
 \le \quad \Delta t \quad \le \quad \frac{1}{6 (1-\theta)} \min\limits_{i}
\frac{|\omega_i|}{\sum\limits_{K\in \omega_i} |K|\, \lambda_{max}(\DDK)\, (h_{i_K}^K)^{-2} },
\label{dt-5}
\eey
where the maximum is taken over all interior edges, $K$ and $K'$ are the two elements sharing the common edge $e_{ij}$,
and the minimum is taken over all interior vertices and $\omega_i$ is the patch of the elements containing
$\V{a}_i$ as its vertex.
\end{thm}

\begin{thm}
\label{thm3.4}
In two dimensions, scheme (\ref{fem-4}) with a lumped mass matrix
satisfies a discrete maximum principle if the mesh satisfies
the Delaunay-type mesh condition (\ref{delaunay-1}) and the time step size satisfies
\bey
\Delta t \quad \le \quad \frac{1}{3 (1-\theta)} \min\limits_{i}
\frac{|\omega_i|}{\sum\limits_{K\in \omega_i} |K|\, \lambda_{max}(\DDK)\, (h_{i_K}^K)^{-2} } .
\label{dt-6}
\eey
\end{thm}

{\Redink{
\begin{rem}
\label{rem3.5}
Conditions (\ref{dt-5}) and (\ref{dt-6}) (for $d=2$) reduce to (\ref{dt-3+5}) and (\ref{dt-4+1}), respectively
for a uniform mesh in the metric specified by $\mathbb{D}^{-1}$ but
are weaker than conditions (\ref{dt-3}) and (\ref{dt-4}) for general meshes.
\proofend
\end{rem}
}}

\section{Numerical results}
\label{Sec-results}

\begin{figure}[thb]
\centering
\hbox{
\begin{minipage}[b]{3in}
\includegraphics[width=3in]{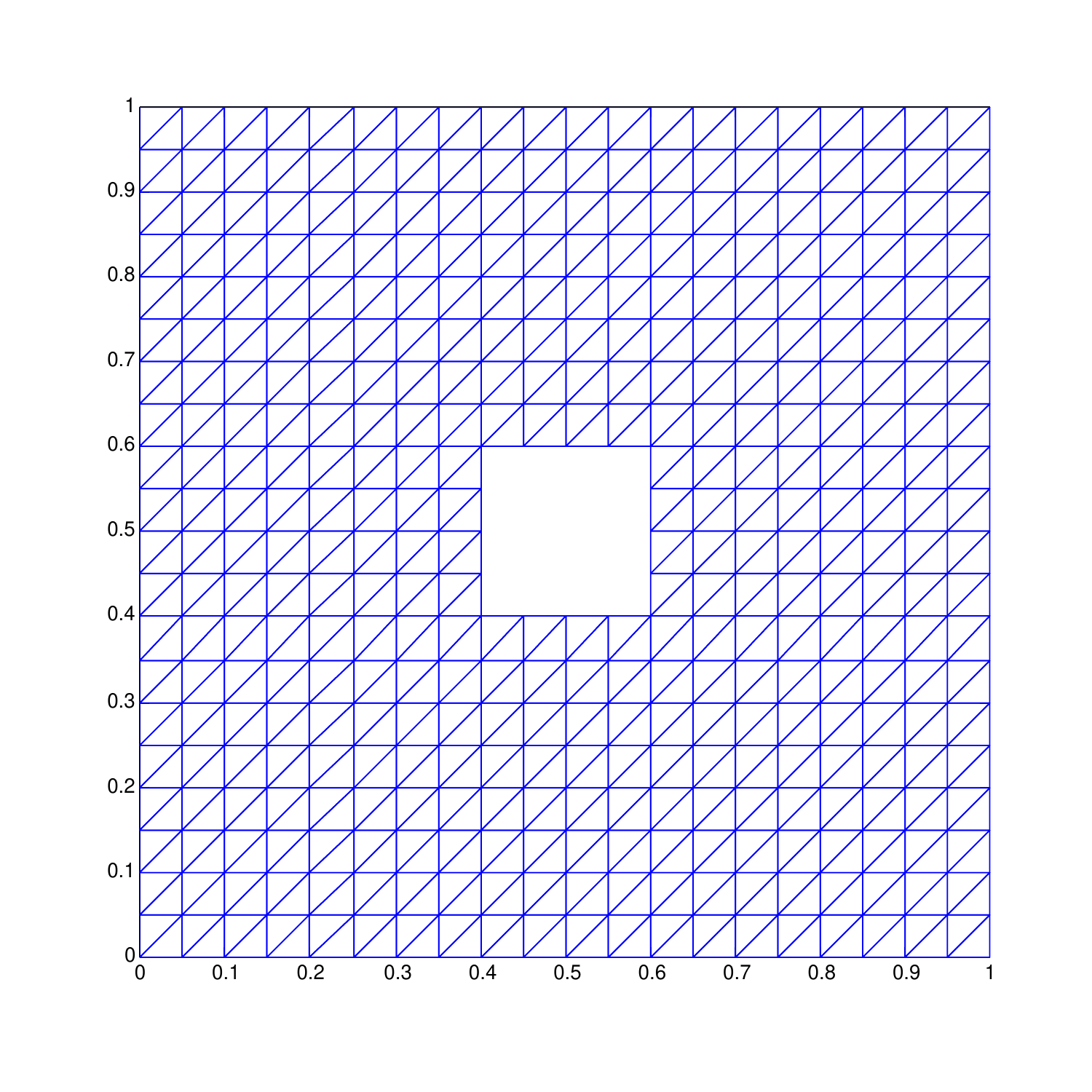}
\centerline{(a): Mesh45}
\end{minipage}
\hspace{10mm}
\begin{minipage}[b]{3in}
\includegraphics[width=3in]{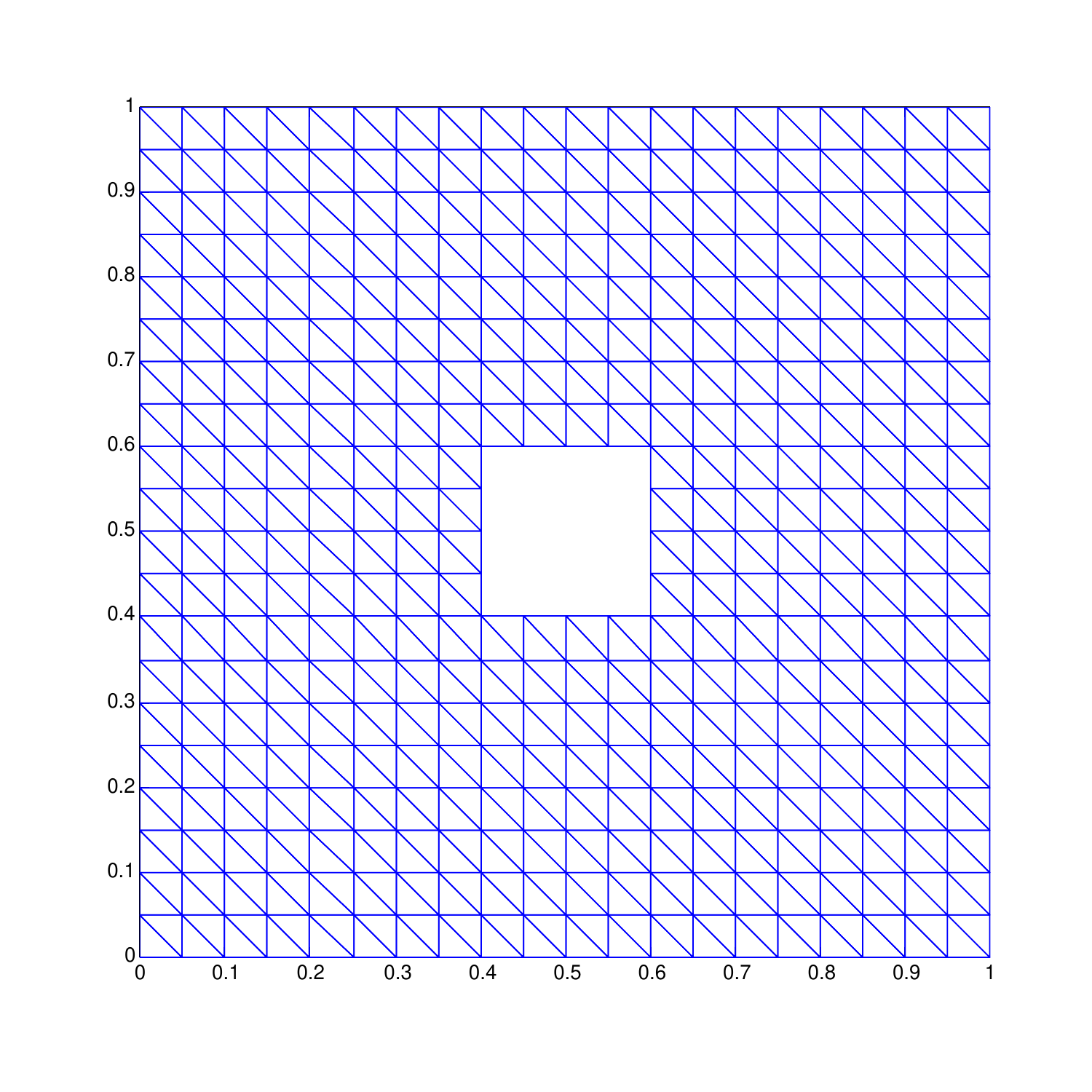}
\centerline{(b): Mesh135}
\end{minipage}
}
\caption{Examples of Mesh45 and Mesh135.}
\label{fixed-mesh}
\end{figure}

In this section we present numerical results obtained for {\Redink{ three }} examples in two dimensions
to demonstrate the significance of both mesh conditions (\ref{anoac-1}) and (\ref{delaunay-1}) and time step conditions (\ref{dt-5}) and (\ref{dt-6}) for DMP satisfaction. Three types of mesh are considered. The first is denoted by Mesh45 where the elements are isosceles right triangles with longest sides in the northeast direction. The second one is denoted by Mesh135 where the elements are isosceles right triangles with longest sides in the northwest direction. Examples of Mesh45 and Mesh135 are shown in Figs. \ref{fixed-mesh}(a) and (b). 
The third type of mesh, denoted by $M_{DMP}$, is a uniform mesh in the metric $M_{DMP} = \mathbb{D}^{-1}$
which guarantees satisfaction of mesh condition (\ref{anoac-1}) (cf. Remark~\ref{rem3.3}).

The implicit Euler method (corresponding to $\theta = 1$ in (\ref{fem-2})) is used in our computation.
For this method, conditions (\ref{dt-3}), (\ref{dt-4}), (\ref{dt-5}), and (\ref{dt-6}) place no constraint on the upper bound of $\Delta t$.
For this reason, we consider only the lower bound for the time step size.
The lower bound in (\ref{dt-3}) (related to the anisotropic nonobtuse angle condition)
is denoted by $\Delta t_{Ani}$ and that in (\ref{dt-5}) (related to the Delaunay-type mesh condition) by $\Delta t_{Del}$.
Unless stated otherwise, the presented results are obtained after 10 steps of time integration.

\begin{exam}
\label{ex1}
The first example is in the form of IBVP (\ref{ibvp-1}) with
\[
f\equiv 0,\quad \Omega = [0,1]^2\backslash [0.4,0.6]^2,
\quad g = 0 \mbox{ on } \Gamma_{out}, \quad  g = 4 \mbox{ on } \Gamma_{in},
\]
where $\Gamma_{out}$ and $\Gamma_{in}$ are the outer and inner boundaries of $\Omega$, respectively;
see Fig. \ref{ex1-domain}(a). The initial solution $u_0(x,y)$ is given as
\[
u_0(x,y) =  \left \{ \begin{array}{ll}
4, & \mbox{ on } \Gamma_{in} \\
0, & \mbox{ in } \Omega \backslash [0.2,0.8]^2  \\
\text{increases linearly}, &  \text{ from } \Gamma_{mid} \text{ to }\Gamma_{in}\\
\end{array} \right.
\]
where $\Gamma_{mid}$ is the boundary of subdomain $[0.2,0.8]^2$; see Fig. \ref{ex1-domain}.
The diffusion matrix is taken as 
\[
\mathbb{D} = \left [\begin{array}{cc} 50.5 & 49.5 \\ 49.5 & 50.5 \end{array} \right ],
\]
which has eigenvalues $100$ and $1$. The principal eigenvectors are in the northeast direction.

This example satisfies the maximum principle and the exact solution (whose analytical expression is unavailable) stays between 0 and 4. Our goal is to produce a numerical solution which also satisfies DMP and stays between 0 and 4. 

\begin{figure}[thb]
\hspace{10mm}
\begin{minipage}[b]{3in}
\begin{center}	
\begin{tikzpicture}[scale = 1.3]
\draw [thick] (0,0) -- (0,4) -- (4, 4) -- (4, 0) -- (0, 0);
\draw [thick] (1.6,1.6) -- (1.6,2.4) -- (2.4, 2.4) -- (2.4, 1.6) -- (1.6, 1.6);
\draw [dashed] (0.8,0.8) -- (3.2,0.8) -- (3.2, 3.2) -- (0.8, 3.2) -- (0.8, 0.8);
\draw (2, -0.35) node {$ u = 0 $};
\draw (4.5, 2) node {$ \Gamma_{out} $};
\draw (2, 1.3) node {$ u = 4 $};
\draw (2.8, 2) node {$ \Gamma_{in} $};
\draw (2, 3.5) node {$ \Gamma_{mid} $};
\end{tikzpicture}
\centerline{(a): Physical domain and boundary}
\end{center}
\end{minipage}
\begin{minipage}[b]{3in}
	\begin{center}
\includegraphics[width=3in]{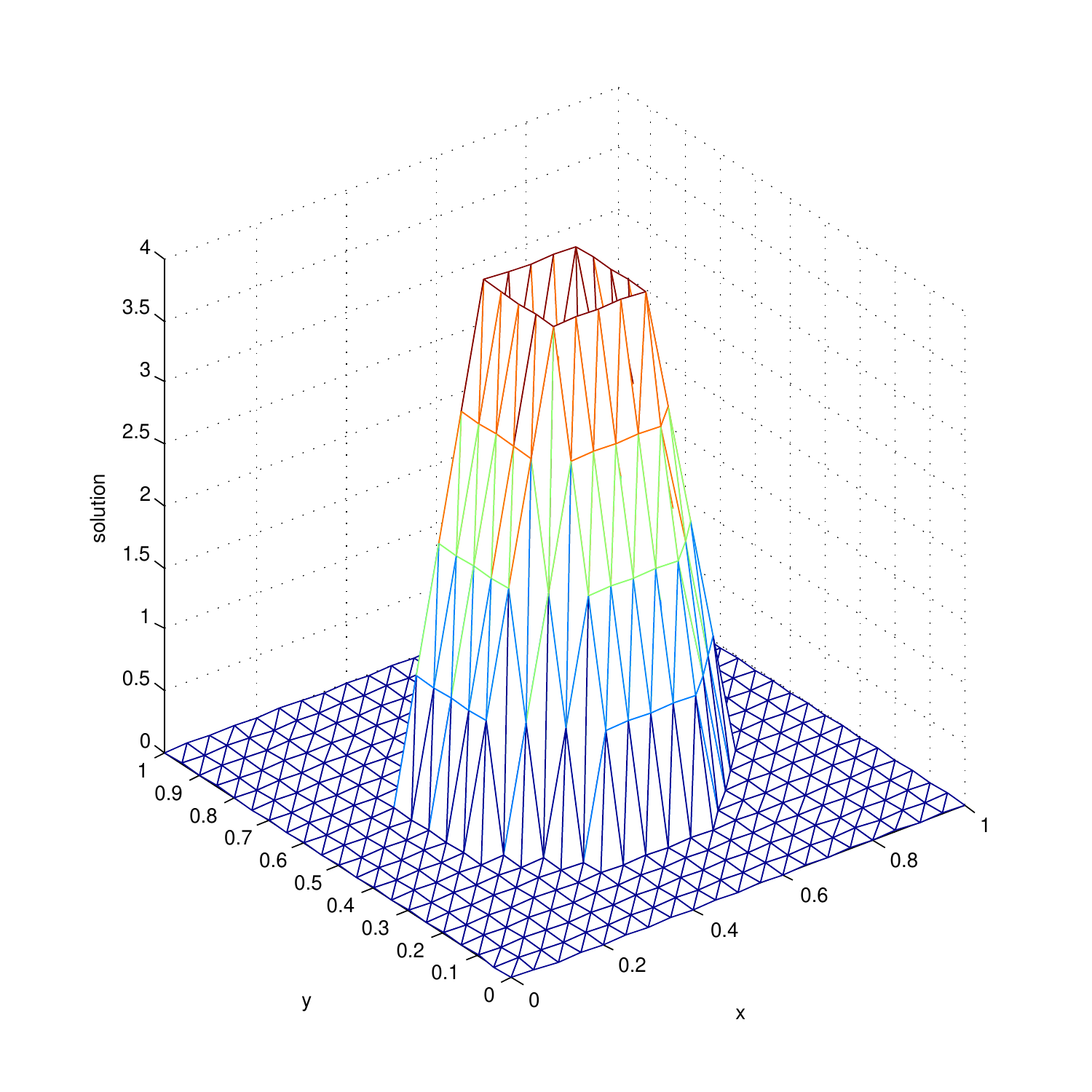}
\centerline{(b): Initial solution}
\end{center}
\end{minipage}
\caption{The physical domain, boundary condition, and initial solution for Example~\ref{ex1}.}
\label{ex1-domain}
\end{figure}

We first consider Mesh45 and Mesh135. Mesh45 satisfies the anisotropic nonobtuse angle condition (\ref{anoac-1}) since
its maximum angle in the metric $M = \DD^{-1}$ is $0.47\pi$. It is known \cite{Hua11} that (\ref{anoac-1}) implies
the Delaunay-type mesh condition, (\ref{delaunay-1}). By direct calculation we can find that the maximum of
the left-hand-side term of  (\ref{delaunay-1}) is $0.94\pi$.
On the other hand, Mesh135 satisfies neither of (\ref{anoac-1}) and (\ref{delaunay-1}),
with the maximum angle in the metric $M = \DD^{-1}$ being $0.94\pi$ and the maximum of the left-hand-side term
of (\ref{delaunay-1}) being $1.87\pi$. 

The solution contours (after 10 time steps) using Mesh45 and Mesh135 with $h=2.5\times 10^{-2}$ and $\Delta t=1.5\times 10^{-4}$ are shown
in Fig. \ref{ex1-soln-fixed}, where $h$ denotes the maximal height of triangular elements of the mesh and $u_{min}$ is the minimum of the numerical solution.  No undershoot occurs in the numerical solution obtained with Mesh45.   

\begin{figure}[thb]
\centering
\hbox{
\begin{minipage}[b]{3in}
\includegraphics[width=3in]{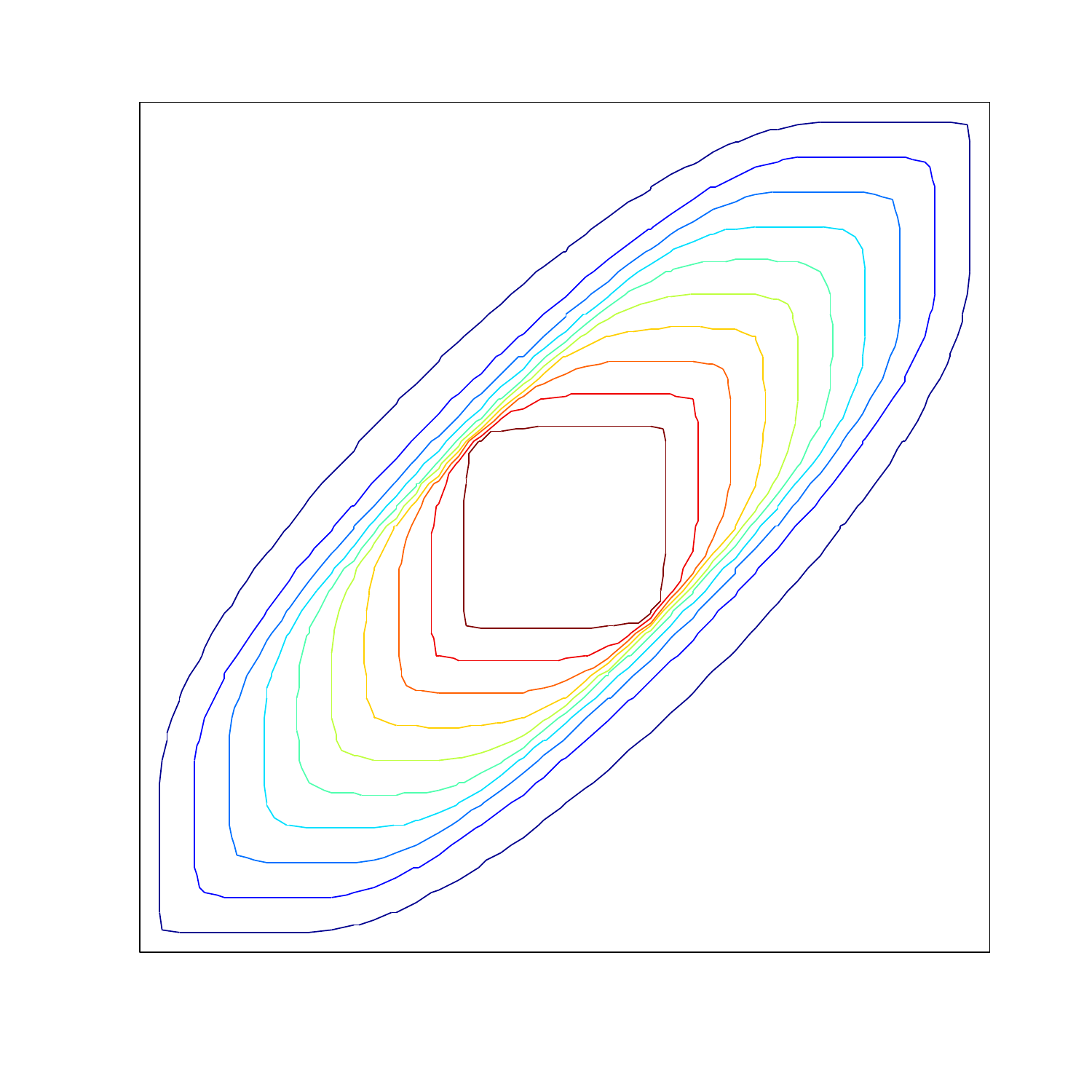}
\centerline{(a): Mesh45, $u_{min}=0$}
\end{minipage}
\hspace{10mm}
\begin{minipage}[b]{3in}
\includegraphics[width=3in]{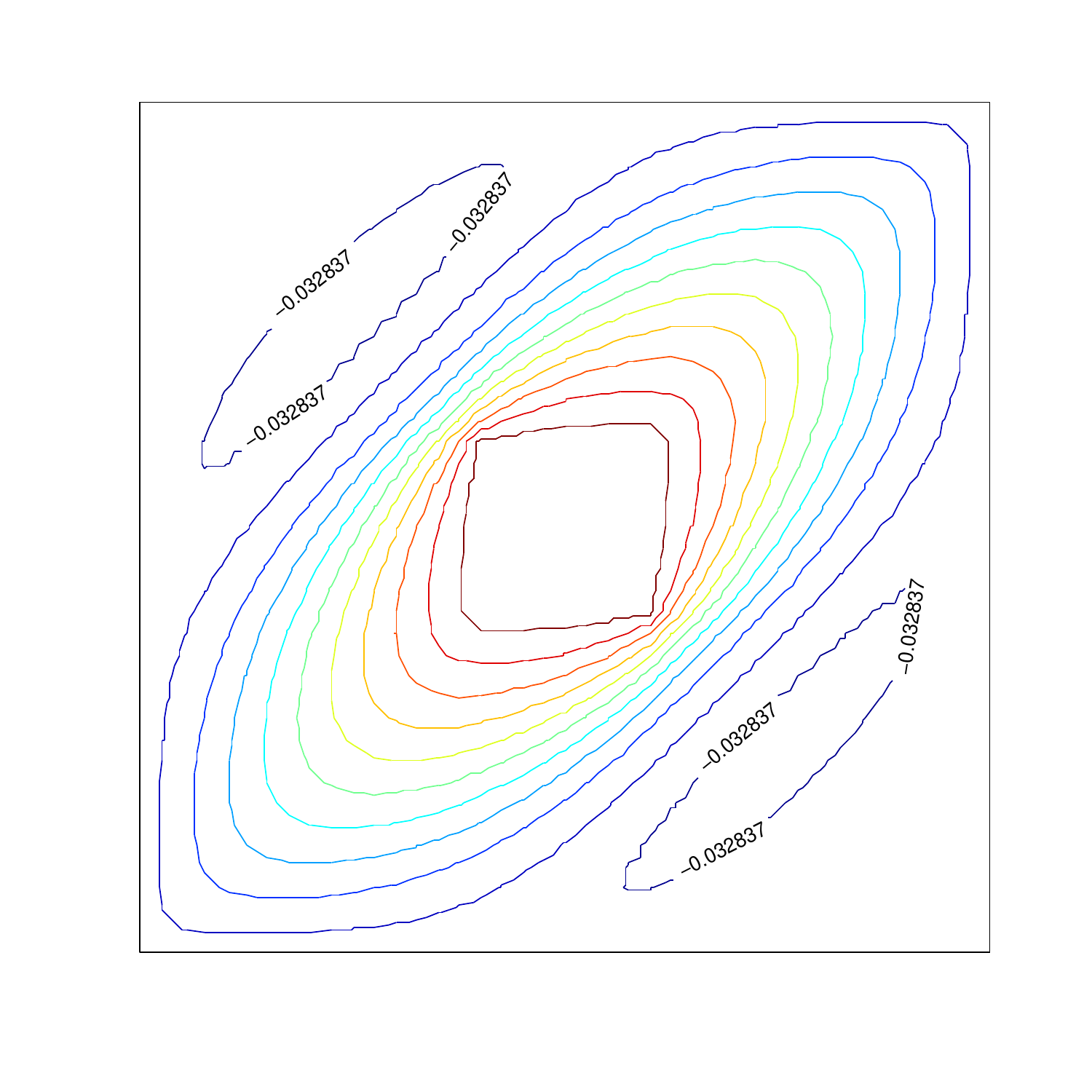}
\centerline{(b): Mesh135, $u_{min}=-6.57\times 10^{-2}$}
\end{minipage}
}
\caption{Solution contours obtained for Mesh45 and Mesh135 with $h=2.5\times 10^{-2}$
and $\Delta t = 1.5\times 10^{-4}$ for Example~\ref{ex1}.}
\label{ex1-soln-fixed}
\end{figure}

The results for Mesh45 are listed in Table \ref{ex1-tab-mesh45}. They show that for meshes with
$h \le 2.5\times 10^{-2}$, $\Delta t_{Del}$ is smaller than the step size $\Delta t=1.5\times 10^{-4}$
used in the computation. As a consequence, time condition (\ref{dt-5}) (and mesh condition (\ref{delaunay-1}))
is satisfied and Theorem \ref{thm3.3}
implies that the numerical solution satisfies DMP. Table~\ref{ex1-tab-mesh45} confirms that
no undershoot occurs in the numerical solution or $u_{min} = 0$.
On the other hand, for $h = 5.0\times 10^{-2}$, neither of time conditions (\ref{dt-3}) and (\ref{dt-5}) is satisfied and undershoot with $u_{min} = -1.41\times 10^{-7}$ is observed.

The table also records the numerical results obtained for $h = 2.5\times 10^{-2}$ and $h = 1.25\times 10^{-2}$ with
decreasing $\Delta t$. One can see that no undershoot occurs when $\Delta t \ge \Delta t_{Del}$. However,
undershoot occurs when $\Delta t$ continues to decrease and pass $\Delta t_{Del}$. This is consistent with
Theorem~\ref{thm3.3}.

It is pointed out that $\Delta t_{Del} < \Delta t_{Ani}$ for all the cases listed in the table.
Moreover, for some cases we have $\Delta t_{Del} < \Delta t < \Delta t_{Ani}$ and no undershoot occurs
in the numerical solution. These indicate that 
time condition (\ref{dt-5}) (related to the Delaunay-type mesh condition) is weaker than
(\ref{dt-3}) (related to the anisotropic nonobtuse angle condition).

Recall that Mesh135 does not satisfy mesh condition (\ref{anoac-1}) nor (\ref{delaunay-1}).
Thus, there is no guarantee that the numerical solution obtained with Mesh135 satisfies DMP.
Indeed, Table~\ref{ex1-tab-mesh135} shows that undershoot occurs in all numerical solutions obtained
with various sizes of Mesh135 and various $\Delta t$.

Next we consider $M_{DMP}$ meshes which are generated as (quasi-)uniform ones in the metric specified by
$M=\DD^{-1}$.  Recall from Remark~\ref{rem3.3} that such meshes satisfy the anisotropic nonobtuse angle
condition (\ref{anoac-1}). In our computation, $M_{DMP}$ meshes are generated using BAMG (bidimensional
anisotropic mesh generator) code developed by Hecht \cite{Hec97}. An example is shown in Fig. \ref{ex1-dmp}(b).
Notice that the elements are aligned with the principal diffusion direction (northeast).
Since the diffusion tensor $\DD$ is constant, the mesh is generated initially based on $M_{DMP} = \DD^{-1}$
and then kept for the subsequent time steps.

The results obtained with $M_{DMP}$ meshes are similar to those obtained with Mesh45. For example,
for the $M_{DMP}$ mesh shown in Fig. \ref{ex1-dmp} (b), it is found numerically that
$\Delta t_{Ani}=4.30\times 10^{-2}$ and $\Delta t_{Del}=1.63\times 10^{-3}$. Theorem~\ref{thm3.3} ensures that
no undershoot occurs in the numerical solution when $\Delta t \ge \Delta t_{Del}$. 
It is emphasized that (\ref{delaunay-1}) and (\ref{dt-5})
are not necessary for DMP satisfaction and the numerical solution may be free of undershoot
for some smaller values of $\Delta t$. In fact, no undershoot is observed numerically for $\Delta t\ge 10^{-4}$.
An undershoot-free solution obtained with the mesh shown in Fig. \ref{ex1-dmp} (b) and time step size
$\Delta t = 1.5\times 10^{-4}$ is shown in Fig. \ref{ex1-dmp} (a). For the same mesh with $\Delta t=1.0\times 10^{-5}$, undershoot is observed with $u_{min}=-1.45\times 10^{-6}$.

\begin{figure}[thb]
\centering
\hbox{
\begin{minipage}[t]{3in}
\includegraphics[width=3in]{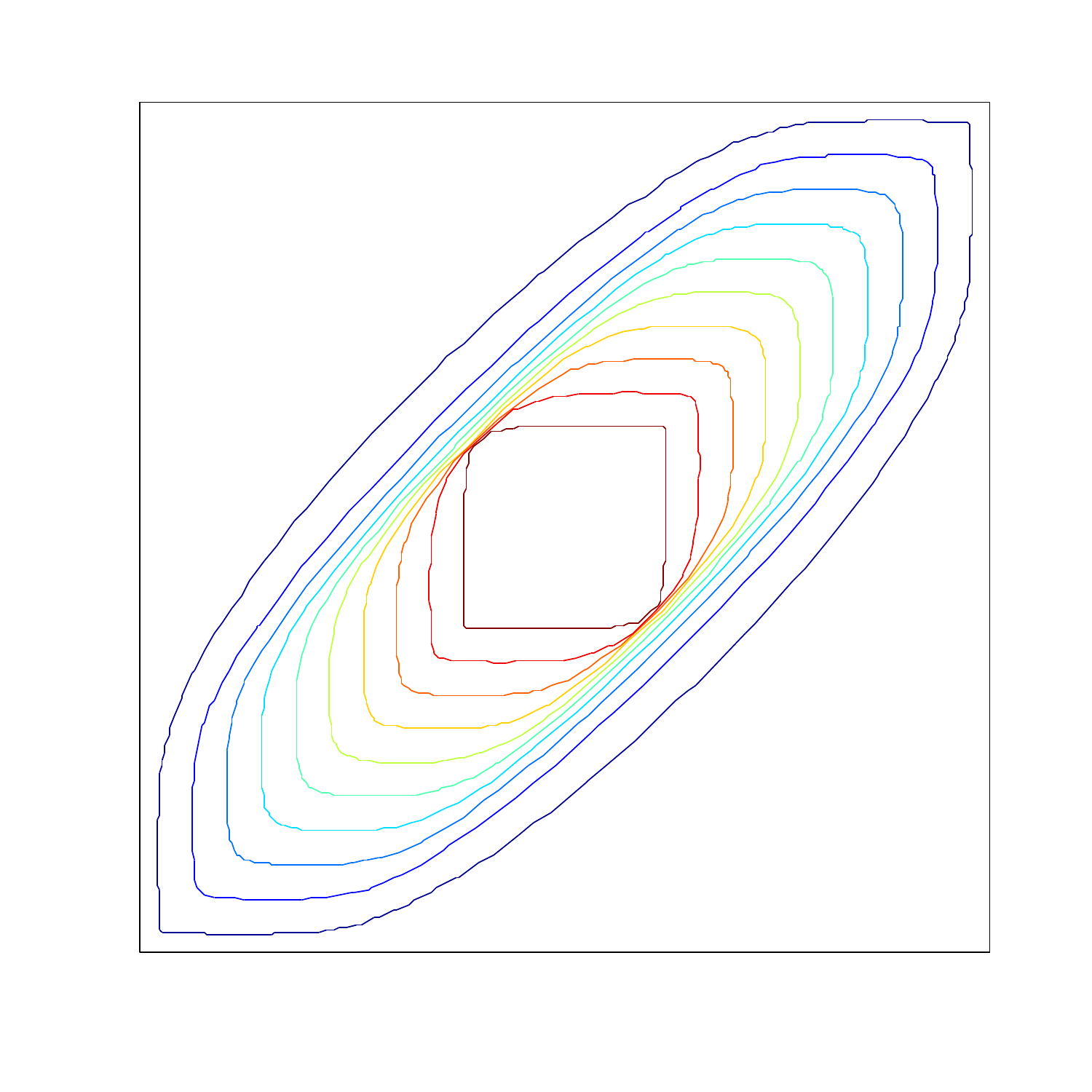}
\centerline{(a): $u_{min}=0$}
\end{minipage}
\hspace{10mm}
\begin{minipage}[t]{3in}
\includegraphics[width=3in]{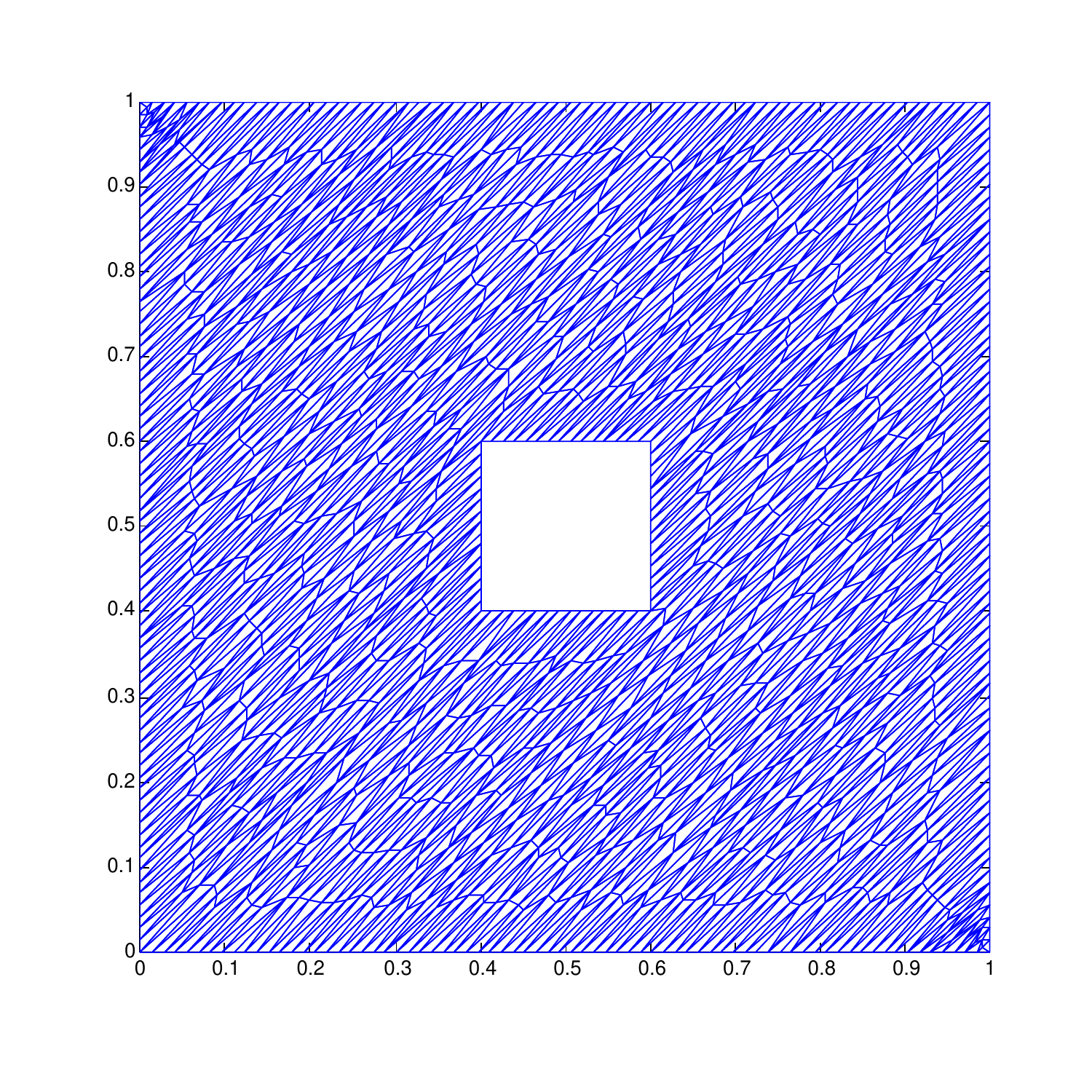}
\centerline{(b): $M_{DMP}$ mesh, $N_e=2362$}
\end{minipage}
}
\caption{An $M_{DMP}$ mesh (with $N_e=2362$ and $N_v=1357$) and
the corresponding solution obtained with $\Delta t = 1.5\times 10^{-4}$ for Example~\ref{ex1}.}
\label{ex1-dmp}
\end{figure}

Finally, we consider the lumped mass method. Theorem~\ref{thm3.4} implies that there is no constraint placed on
$\Delta t$ for the DMP satisfaction of the numerical solution with the lumped mass matrix and implicit Euler
discretization. Indeed, for all Mesh45 meshes and $\Delta t$ considered in Table \ref{ex1-tab-mesh45},
no undershoot is observed numerically for the lumped mass method. The same also holds for $M_{DMP}$ meshes.
For example, for the mesh shown in Fig. \ref{ex1-dmp}(b), no undershoot is observed in the numerical solution
for $\Delta t=10^{-4}$, $10^{-5}$, and $10^{-6}$.
For Mesh135 meshes, mesh condition (\ref{anoac-1}) or (\ref{delaunay-1}) is not satisfied and thus Theorem~\ref{thm3.4}
does not hold. For example, for a case with a Mesh135 mesh with $h=1.25\times 10^{-2}$ and $\Delta t=1.5\times 10^{-4}$, the numerical solution violates DMP and has a minimum $u_{min}=-1.60\times 10^{-2}$. 
\end{exam}

\begin{exam}
\label{ex2}
The second example is the same as Example \ref{ex1} except that the diffusion matrix is
taken as a function of $x$ and $y$, i.e., 
\begin{equation}
\DD = \left (\begin{array}{cc} \cos\theta & -\sin\theta \\ \sin\theta & \cos\theta \end{array} \right )
 \left (\begin{array}{cc} k_1 & 0 \\ 0 & k_2 \end{array} \right )
\left (\begin{array}{cc} \cos\theta & \sin\theta \\ -\sin\theta & \cos\theta \end{array} \right ),
\label{D-1}
\end{equation}
where $k_1= 100$, $k_2 = 1$, and 
$\theta = \theta(x,y)$ is the angle of the tangential direction at point $(x,y)$ along circles centered at
$(0.5,0.5)$. This diffusion matrix $\DD$ also has eigenvalues $1$ and $100$ but has its principal
eigen-direction along the tangential direction of circles centered at $(0.5,0.5)$.
A physical example with such a diffusion matrix is the toroidal magnetic field in a Tokamak device
confining fusion plasma \cite{SCTPKB10}. This problem also satisfies the maximum principle and
the solution stays between 0 and 4. 

For this example, neither Mesh45 nor Mesh135 (cf. Fig. \ref{fixed-mesh}) satisfies the 
Delaunay-type mesh condition (\ref{delaunay-1}). In the metric specified by $M=\DD^{-1}$,
the maximum of the left-hand side of the inequality is $1.87\pi$ for both Mesh45 and Mesh135.
Due to the symmetry of the diffusion matrix, both Mesh45 and Mesh135 lead to almost the same results
for this example except that undershoot occurs at different locations. Fig. \ref{ex2-fig} shows the results
obtained with these meshes for $\Delta t=5\times 10^{-5}$.

\begin{figure}[thb]
\centering
\hbox{
\hspace{10mm}
\begin{minipage}[thb]{2.5in}
\includegraphics[width=2.5in]{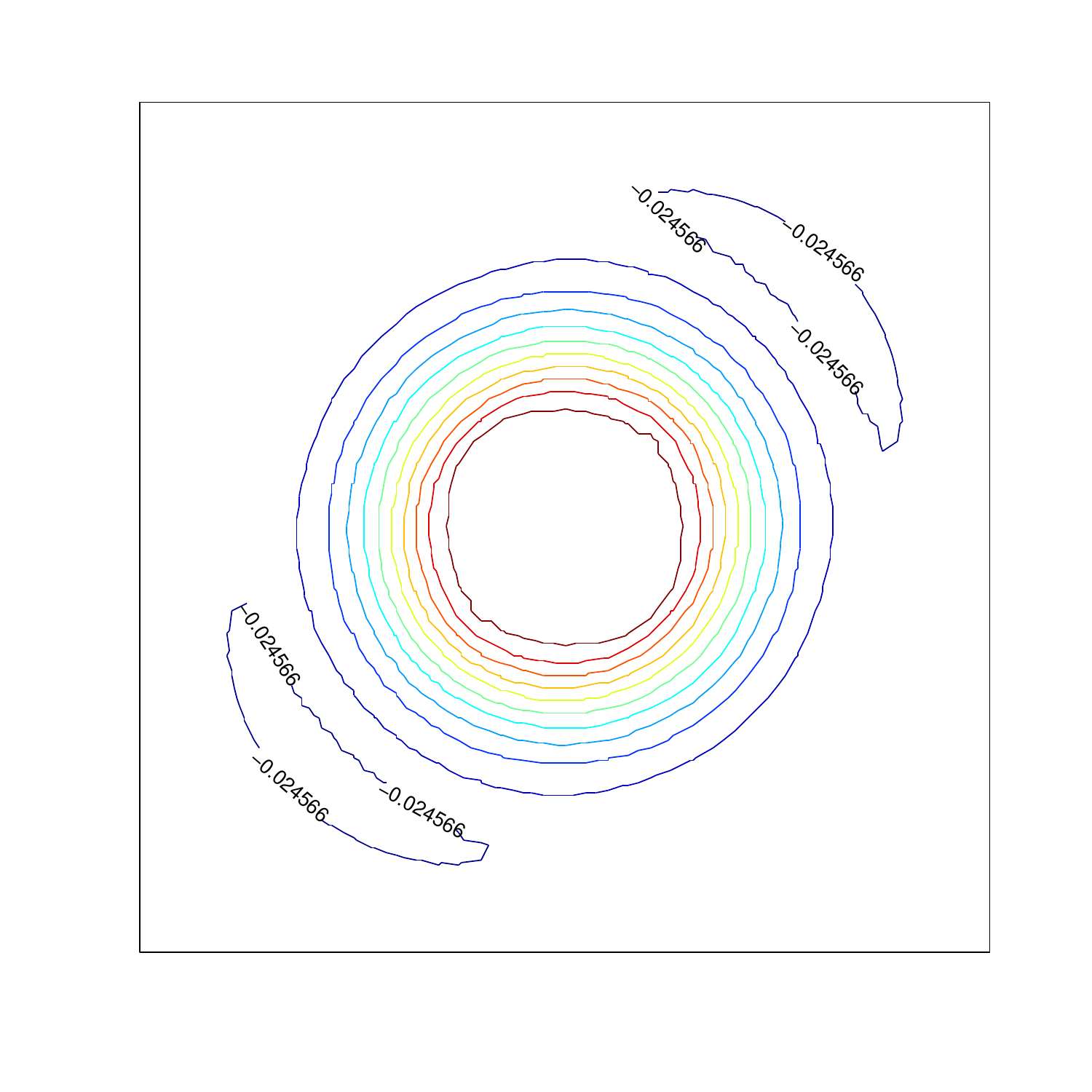}
\centerline{(a): Mesh45, $u_{min}=-4.91\times 10^{-2}$}
\end{minipage}
\hspace{10mm}
\begin{minipage}[thb]{2.5in}
\includegraphics[width=2.5in]{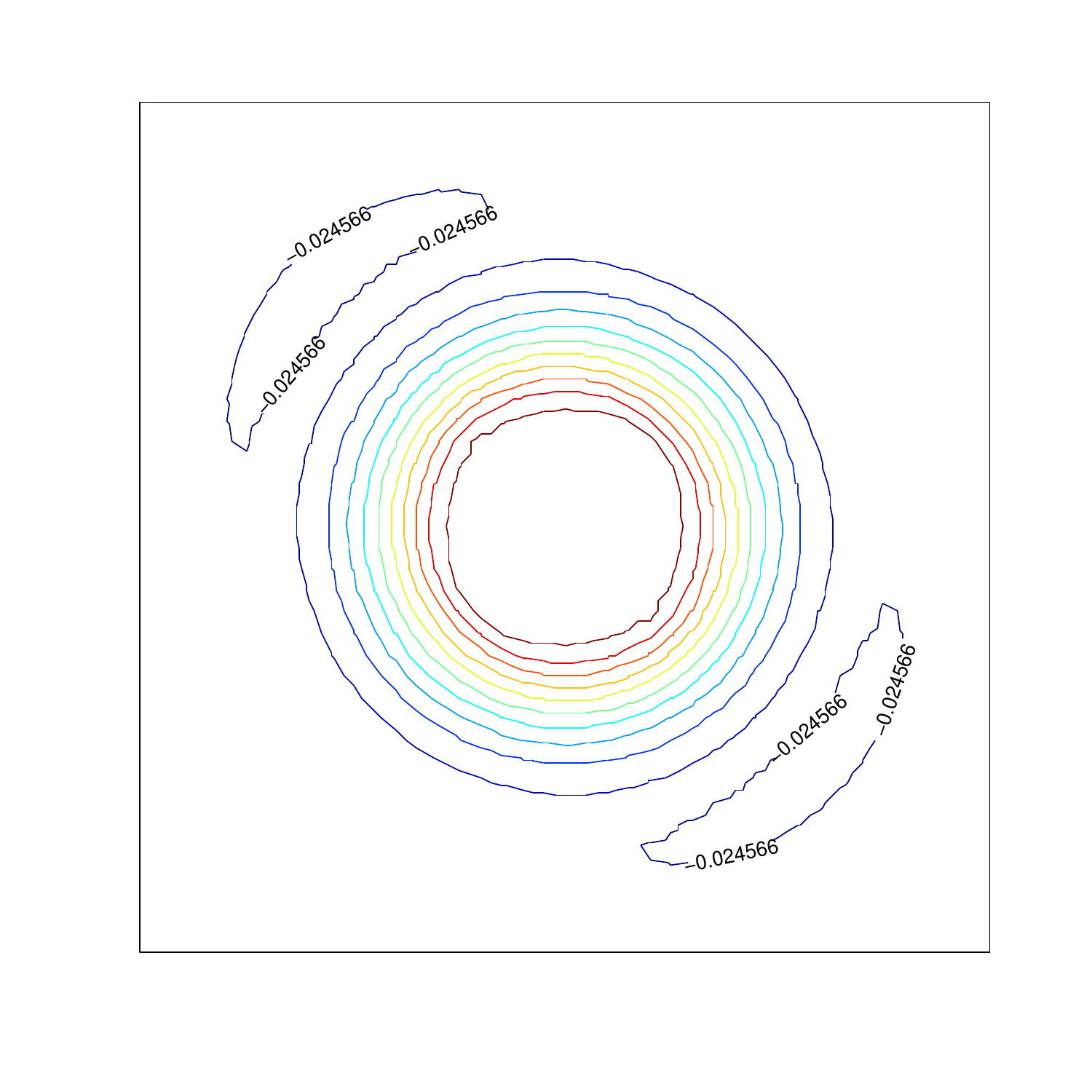}
\centerline{(b): Mesh135, $u_{min}=-4.91\times 10^{-2}$}
\end{minipage}
}
\hbox{
\hspace{10mm}
\begin{minipage}[thb]{2.5in}
\includegraphics[width=2.5in]{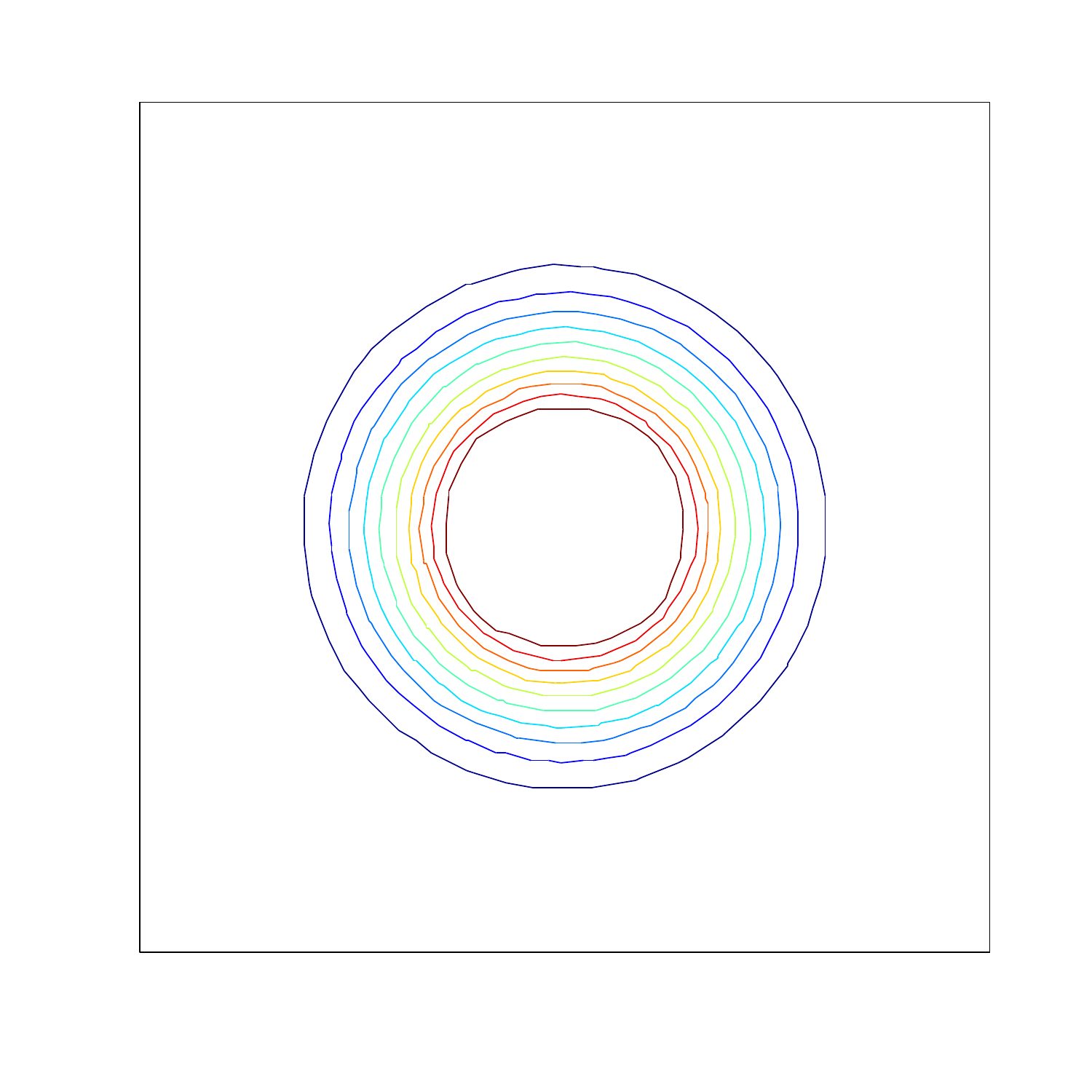}
\centerline{(c):  $M_{DMP}$, $u_{min}=0$}
\end{minipage}
\hspace{10mm}
\begin{minipage}[thb]{2.5in}
\includegraphics[width=2.5in]{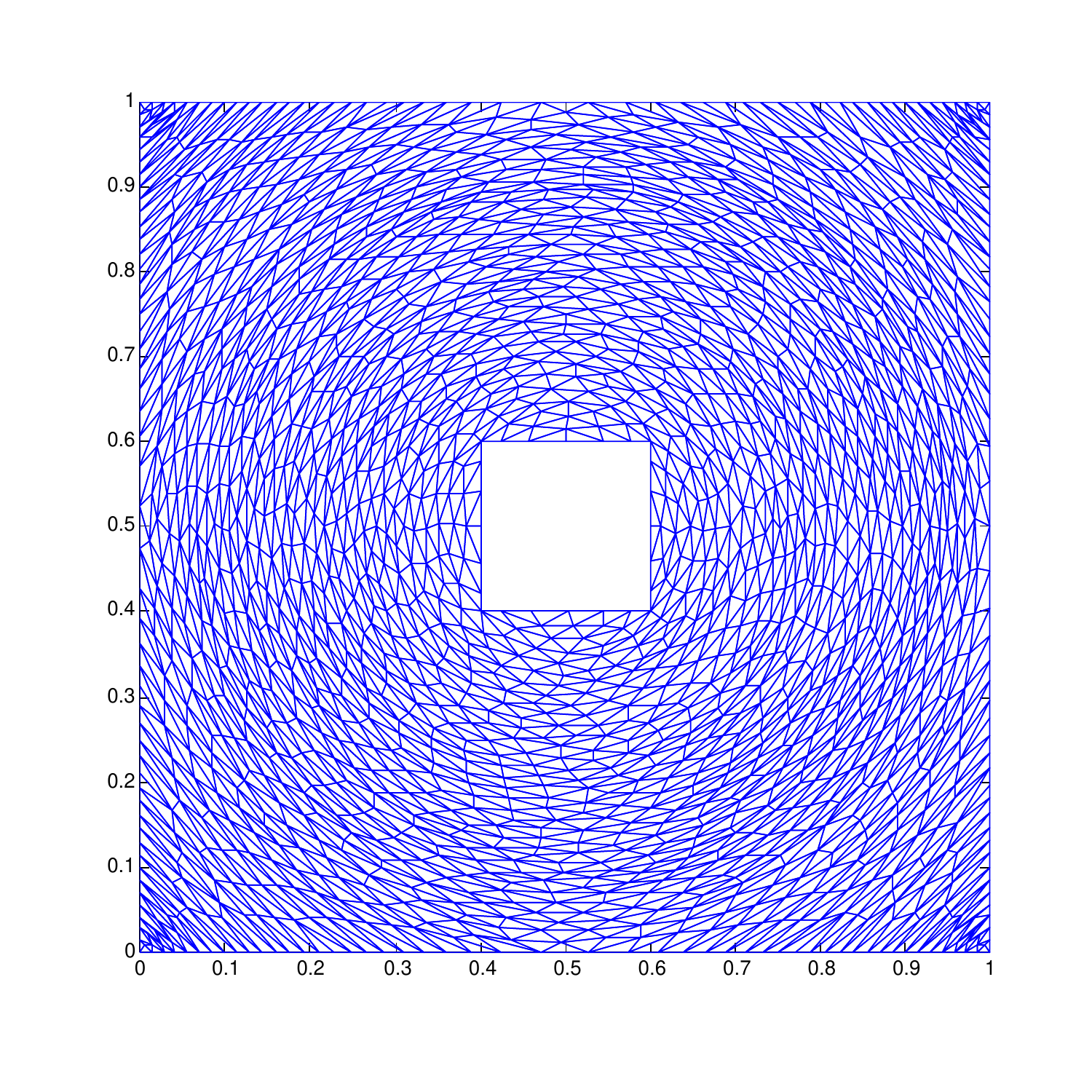}
\centerline{(d): $M_{DMP}$, $N_e=3381$}
\end{minipage}
}
\caption{Results obtained with $\Delta t = 5\times 10^{-5}$ for Example~\ref{ex2}.}
\label{ex2-fig}
\end{figure}

Table \ref{ex2-tab} lists numerical results obtained with Mesh45 and $M_{DMP}$ meshes. 
Recall that Theorem~\ref{thm3.3} does not apply to Mesh45 meshes since they do not satisfy (\ref{delaunay-1}).
As a matter of fact, numerical solutions obtained with this type of meshes with or without mass lumping violate
DMP and exhibit undershoot. On the other hand, $M_{DMP}$ meshes generated with $M=\DD^{-1}$
satisfy the mesh condition. For the lumped mass method, no undershoot occurs in the numerical solution for
all values of $\Delta t$. This is consistent with Theorem~\ref{thm3.4}. For the standard finite element method,
there is no undershoot for relatively large $\Delta t$. It is interesting to point out that for this example with
variable $\DD$, the lower bounds $\Delta t_{Ani}$ and $\Delta t_{Del}$ are far too pessimistic. A several magnitude smaller $\Delta t$ can still lead to numerical solutions free of undershoot.

\end{exam}

{\Redink{
\begin{exam}
\label{ex3}
This example is the same as the previous examples except that the diffusion matrix is
taken as in the form (\ref{D-1}) with
\[
\theta = \frac{1}{2}\arctan(\cos(\frac{\pi x}{4})),\quad
k_1 = 100 \cos((x^2+y^2)\frac{\pi}{6}),\quad
k_2 = 10 \sin((x^2+y^2+1)\frac{\pi}{6}).
\]
Notice that $\DD$ is a function of $x$ and $y$ and
both its eigenvalues and eigenvectors vary with location.

Numerical results are shown in Table \ref{ex3-tab} and
Fig. \ref{ex3-fig}. Similar observations can be made as in the previous example. More specifically,
both Mesh45 and Mesh135 does not satisfy the Delaunay-type mesh condition (\ref{delaunay-1})
and thus there is no guarantee that the obtained numerical solution is undershoot-free.
On the other hand, $M_{DMP}$ meshes generated with $M=\DD^{-1}$
satisfy (\ref{delaunay-1}). The numerical solution is guaranteed to be undershoot-free
for sufficiently large $\Delta t$ for the standard linear finite element method and for all $\Delta t$
for the lumped mass method.

\begin{figure}[thb]
\centering
\hbox{
\hspace{10mm}
\begin{minipage}[thb]{2.5in}
\includegraphics[width=2.5in]{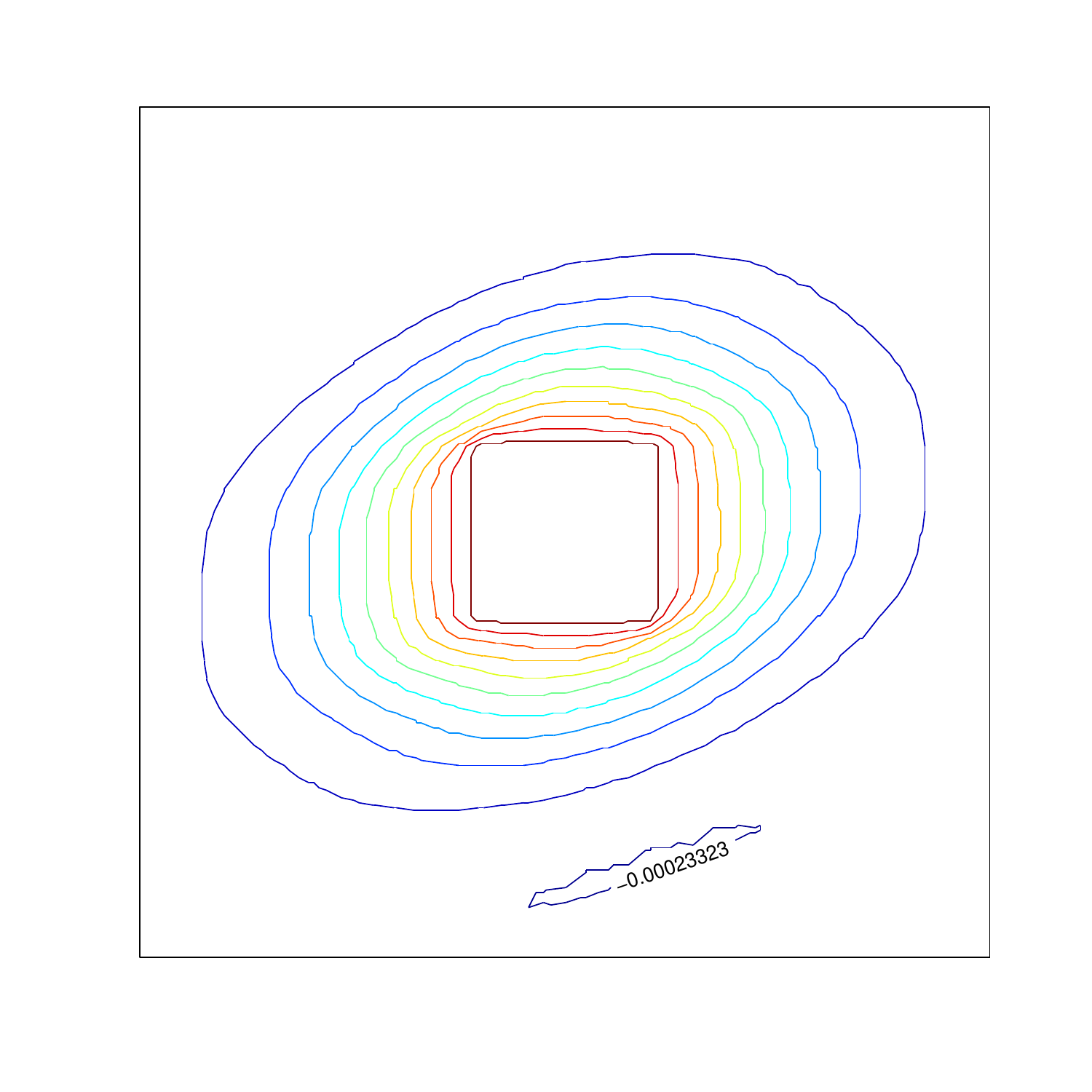}
\centerline{(a): Mesh45, $u_{min}=-4.66 \times 10^{-4}$}
\end{minipage}
\hspace{10mm}
\begin{minipage}[thb]{2.5in}
\includegraphics[width=2.5in]{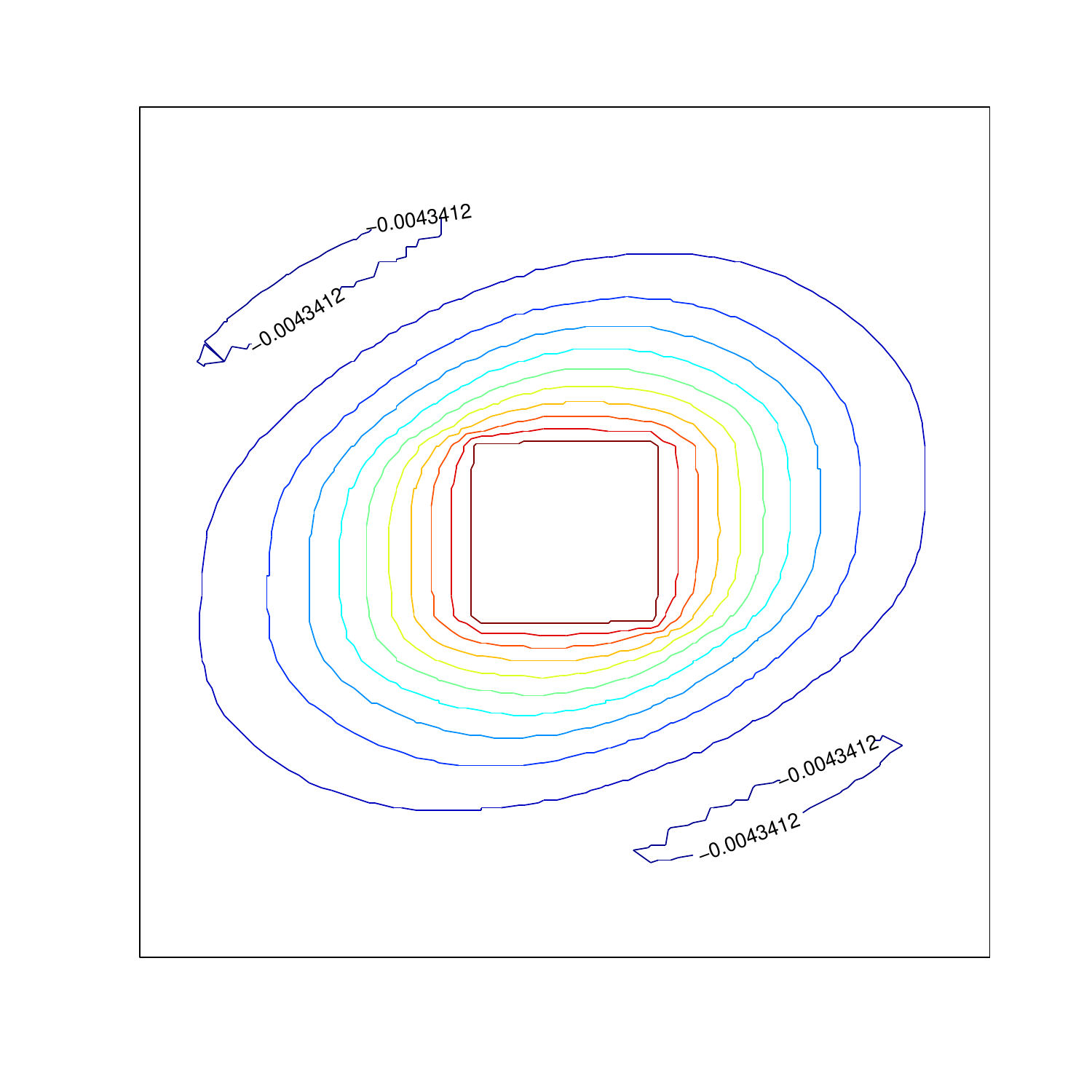}
\centerline{(b): Mesh135, $u_{min}=-8.68 \times 10^{-3}$}
\end{minipage}
}
\hbox{
\hspace{10mm}
\begin{minipage}[thb]{2.5in}
\includegraphics[width=2.5in]{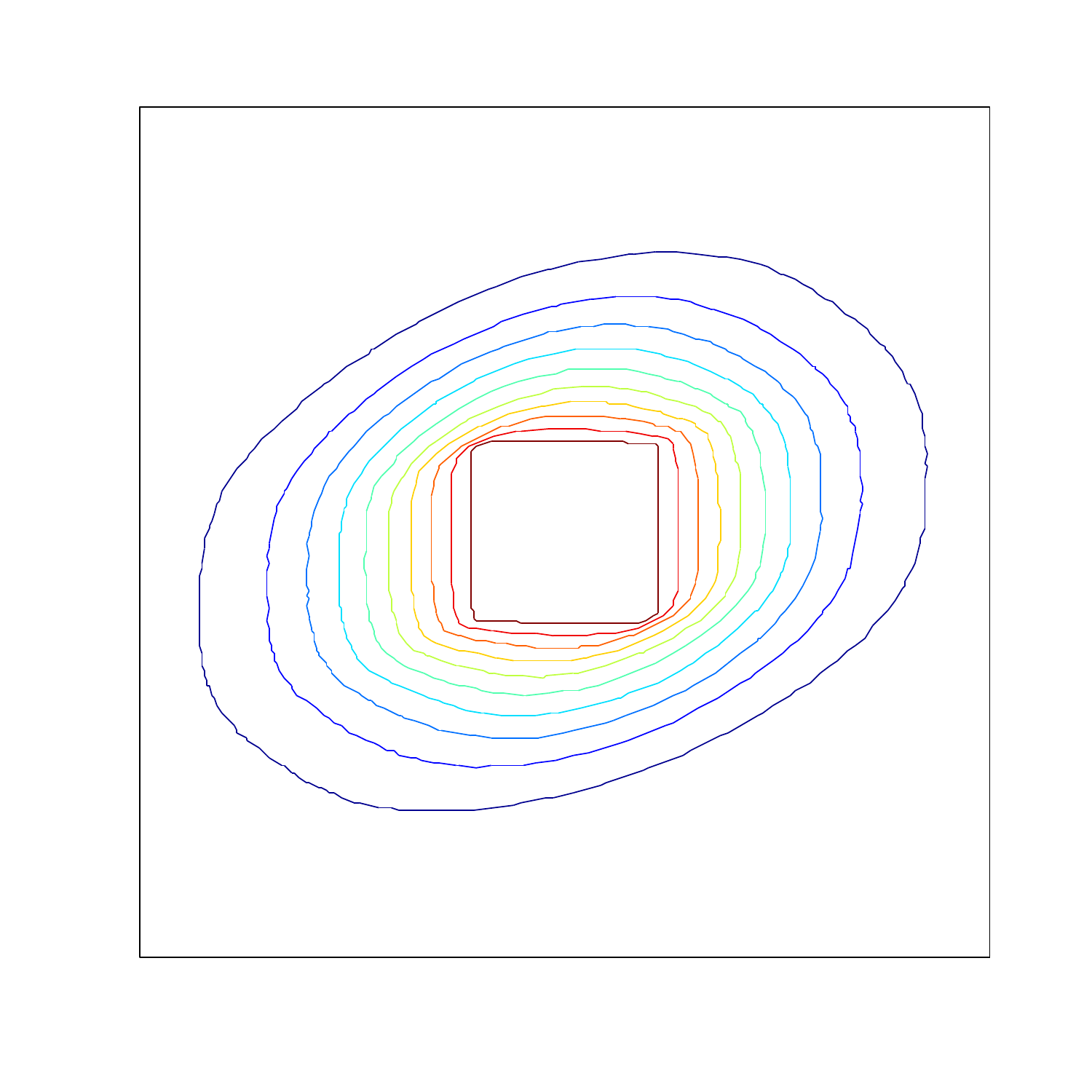}
\centerline{(c):  $M_{DMP}$, $u_{min}=0$}
\end{minipage}
\hspace{10mm}
\begin{minipage}[thb]{2.5in}
\includegraphics[width=2.5in]{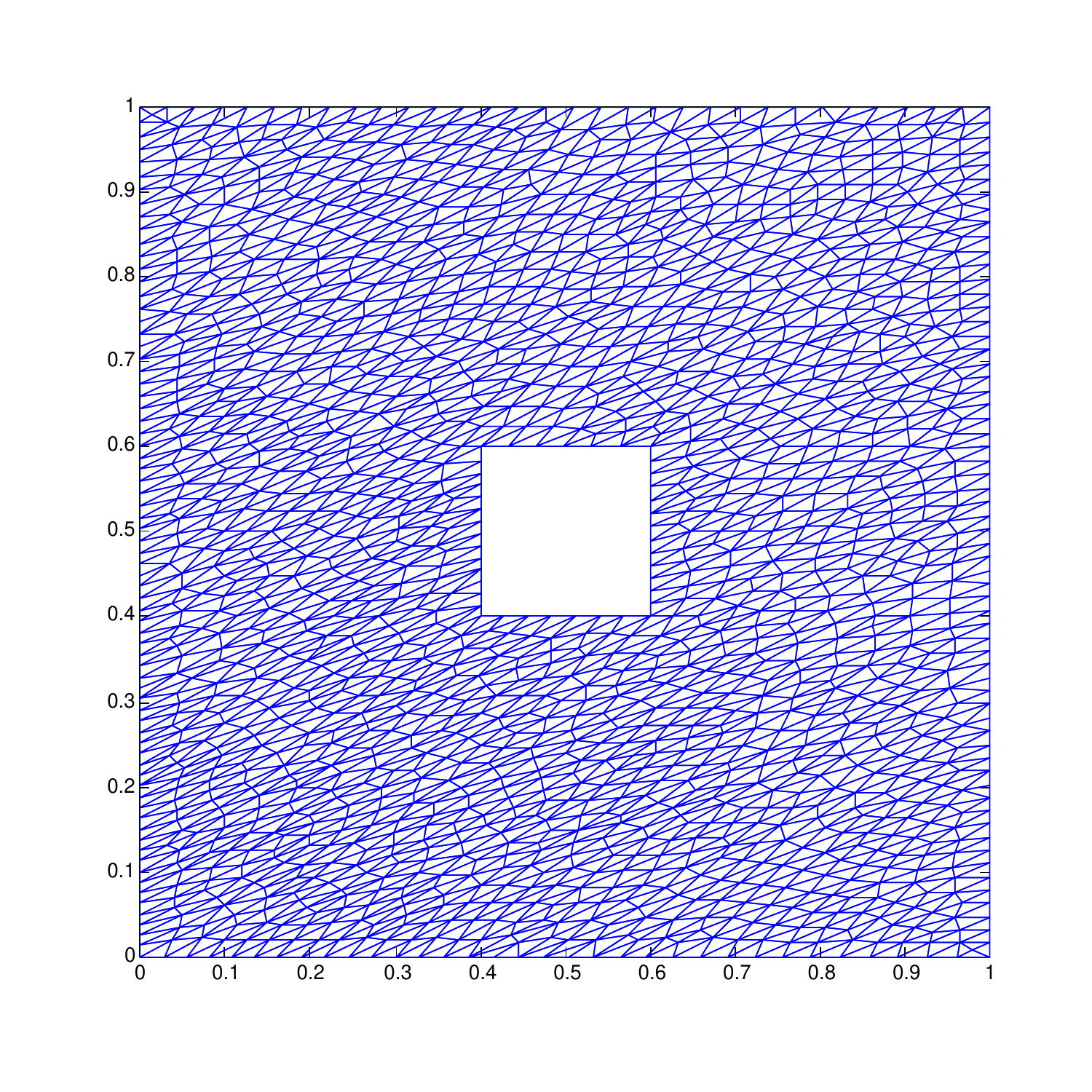}
\centerline{(d): $M_{DMP}$, $N_e=3180$}
\end{minipage}
}
\caption{Results obtained with $\Delta t = 1 \times 10^{-5}$ for Example~\ref{ex3}.}
\label{ex3-fig}
\end{figure}

\end{exam}
}}

\section{Conclusions}
\label{Sec-con}

In the previous sections we have studied the conditions under which a full discretization for IBVP
 (\ref{ibvp-1}) with a general diffusion matrix function
 satisfies a discrete maximum principle.
 The discretization is realized using the $\theta$-method in time and the linear finite element method
 in space. The main theoretical results are given in Theorems \ref{thm3.1}, \ref{thm3.2}, \ref{thm3.3}, and \ref{thm3.4}.
 
Specifically, the numerical solution obtained with the full discrete scheme satisfies a discrete maximum
 principle when the mesh satisfies the anisotropic nonobtuse angle condition (\ref{anoac-1}) and the time
 step size satisfies condition (\ref{dt-3}). As shown in \cite{LH10}, a mesh satisfying (\ref{anoac-1}) can be generated
 as a uniform mesh in the metric specified by $\alpha\; \DD^{-1}$ with $\alpha$ being a scalar
 function defined on $\Omega_T$. On the other hand, condition (\ref{dt-3}) essentially requires
 the time step size to satisfy
 \begin{equation}
 C_1 h^2 \le \Delta t \le \frac{C_2}{1-\theta} h^2,
 \label{dt-10}
 \end{equation}
 where $C_1$ and $C_2$ are positive constants, $h$ is the maximal element diameter, and $\theta \in (0,1]$
 is the parameter used in the $\theta$-method. Obviously, this condition is restrictive. This is especially true
 when the numerical scheme with $\theta \in [0.5, 1]$ is known to be unconditionally stable and no constraint is placed
 on $\Delta t$ for the sake of stability. Moreover, {\Redink{the presence of the lower bound for $\Delta t$
and the numerical results showing the violation of the maximum principle as $\Delta t \to 0$
seem to support}}
the finding of  Thom{\'e}e and Wahlbin \cite{TW08} that a semi-discrete standard Galerkin finite element solution
 violates DMP since the semi-discrete scheme can be considered as the limit of
 the full discrete scheme as $\Delta t \to 0$.
 Furthermore, Theorems \ref{thm3.2} and \ref{thm3.4} show
 that the lower bound requirement on $\Delta t$ can be removed when a lumped
 mass matrix is used. Finally, in two dimensions, the mesh and time step conditions can be replaced
 with weaker conditions (\ref{delaunay-1}) and (\ref{dt-5}), respectively. Numerical results
 in Sect.~\ref{Sec-results} confirm the theoretical findings.

\vspace{20pt}

{\bf Acknowledgment.} This work was supported in part by NSF under grant DMS-1115118.   


\newpage

\noindent {\bf List of Figures}

\begin{itemize}
	\item[] Figure \ref{f-0}: Sketch of coordinate transformations from $\hat K$ to $K$ and to $\widetilde{K}$.
	\item[] Figure \ref{fixed-mesh}: Examples of Mesh45 and Mesh135.
	\item[] Figure \ref{ex1-domain}: The physical domain, boundary condition, and initial solution for Example~\ref{ex1}.
	\item[] Figure \ref{ex1-soln-fixed}: Solution contours obtained for Mesh45 and Mesh135 with $h=2.5\times 10^{-2}$ and $\Delta t = 1.5\times 10^{-4}$ for Example~\ref{ex1}.
	\item[] Figure \ref{ex1-dmp}: An $M_{DMP}$ mesh (with $N_e=2362$ and $N_v=1357$) and
	the corresponding solution obtained with $\Delta t = 1.5\times 10^{-4}$ for Example~\ref{ex1}.
	\item[] Figure \ref{ex2-fig}: Results obtained with $\Delta t = 5\times 10^{-5}$ for Example~\ref{ex2}.
	\item[] Figure \ref{ex3-fig}: Results obtained with $\Delta t = 1 \times 10^{-5}$ for Example~\ref{ex3}.
\end{itemize}

\newpage

\noindent {\bf Tables}

\begin{table}[bht]
\caption{Numerical results obtained with Mesh45 for Example \ref{ex1}.}
\vspace{2pt}
\centering
\begin{tabular}{cccccc}
\hline\hline
$h$ & $\Delta t_{Ani}$ & $\Delta t_{Del}$ & $\Delta t$ & $u_{min}$ \\ [0.5ex] 
\hline 
5.0e-2 & 1.48e-3 & 3.79e-4 & 1.5e-4 & -1.41e-7 \\
2.5e-2 & 3.70e-4 & 9.47e-5 & 1.5e-4 & 0 \\
1.25e-2 & 9.25e-5 & 2.37e-5 & 1.5e-4 & 0 \\
6.25e-3 & 2.31e-5 & 5.92e-6 & 1.5e-4 & 0 \\
3.125e-3 & 5.78e-6 & 1.48e-6 & 1.5e-4 & 0 \\ 
\hline
2.5e-2 & 3.70e-4 & 9.47e-5 & 1.5e-4 & 0 \\
2.5e-2 & 3.70e-4 & 9.47e-5 & 1.0e-4 & 0 \\
2.5e-2 & 3.70e-4 & 9.47e-5 & 5.0e-5 & -7.91e-10 \\
1.25e-2 & 9.25e-5 & 2.37e-5 & 1.5e-4 & 0 \\
1.25e-2 & 9.25e-5 & 2.37e-5 & 1.0e-5 & -1.31e-6 \\ [1ex]
\hline \hline
\end{tabular}
\label{ex1-tab-mesh45}
\end{table}

\begin{table}[bht]
\caption{Numerical results obtained with Mesh135 for Example \ref{ex1}.}
\vspace{2pt}
\centering
\begin{tabular}{cccccc}
\hline\hline
$h$ & $\Delta t_{Ani}$ & $\Delta t_{Del}$ & $\Delta t$ & $u_{min}$ \\ [0.5ex] 
\hline
5.0e-2 & 1.48e-4 & 2.08e-6 & 1.5e-4 & -8.99e-2 \\
2.5e-2 & 3.70e-5 & 5.21e-7 & 1.5e-4 & -6.57e-2 \\
1.25e-2 & 9.25e-6 & 1.30e-7 & 1.5e-4 & -1.58e-2 \\
1.25e-2 & 9.25e-6 & 1.30e-7 & 1.0e-7 & -2.26e-2 \\
6.25e-3 & 2.31e-6 & 3.26e-8 & 5.0e-4 & -1.59e-3 \\
6.25e-3 & 2.31e-6 & 3.26e-8 & 1.5e-5 & -1.43e-2 \\
6.25e-3 & 2.31e-6 & 3.26e-8 & 1.5e-6 & -2.11e-2 \\ [1ex]
\hline \hline
\end{tabular}
\label{ex1-tab-mesh135}
\end{table}

\begin{table}[bht]
\caption{Results obtained with Mesh45 and $M_{DMP}$ meshes for Example \ref{ex2}.}
\vspace{2pt}
\centering
\begin{tabular}{ccccccc}
\hline\hline
Mesh & $N_e$ & $\Delta t_{Ani}$ & $\Delta t_{Del}$ & $\Delta t$ & $u_{min}$ & $u_{min}$ (lumped mass)\\ [0.5ex]
\hline
Mesh45 & 3072 & 3.47e-4 & 1.17e-2 & 1.0e-4 & -4.31e-2 & -4.11e-2 \\
 &  &  &  & 5.0e-5 & -4.91e-2 & -4.78e-2 \\
&  &  &  & 2.0e-5 & -5.49e-2 &-5.36e-2 \\
 &  &  &  & 1.0e-5 & -5.70e-2 & -5.26e-2 \\
\hline
$M_{DMP}$ & 3381 & 8.61e-2 & 3.06e-2 & 5.0e-2 & 0 & 0 \\
&  &  &  &1.0e-4 & 0 & 0 \\
&  &  &  & 5.0e-5 & 0 & 0 \\
&  &  &  &  2.0e-5 & -1.20e-5  & 0 \\ 
 &  &  &  & 1.0e-5 & -5.02e-4 & 0 \\ [1ex]
\hline \hline
\end{tabular}
\label{ex2-tab}
\end{table}

\begin{table}[thb]
\caption{Results obtained with $M_{DMP}$ meshes for Example \ref{ex3}.}
\vspace{2pt}
\centering
\begin{tabular}{cccccc}
\hline\hline
 $N_e$ & $\Delta t_{Ani}$ & $\Delta t_{Del}$ & $\Delta t$ & $u_{min}$ & $u_{min}$ (lumped mass)\\ [0.5ex]
\hline
 3180 & 1.83e-2 & 6.38e-4 & 1.0e-4 & 0 & 0 \\
 &  &  & 5.0e-5 & 0 & 0 \\
 &  &  & 1.0e-5 & 0 & 0 \\
 &  &  & 2.5e-6 & -7.67e-5  & 0 \\ 
 &  &  & 1.0e-6 & -6.21e-3 & 0 \\ [1ex]
\hline \hline
\end{tabular} 
\label{ex3-tab}
\end{table}


\begin{thebibliography}{10}

{\Redink{\bibitem{FPPGW09}
P.~E.~Farrell, M.~D.~Piggott, C.~C.~Pain, G.~J.~Gorman, and C.~R. Wilson.
\newblock Conservative interpolation between unstructured meshes via supermesh construction.
\newblock {\em Comput. Methods Appl. Mech. Engrg.} 198:2632--2642, 2009.}}

\bibitem{GL09}
S.~G\H{u}nter and K.~Lackner.
\newblock A mixed implicit-explicit finite difference scheme for heat transport
  in magnetised plasmas.
\newblock {\em J. Comput. Phys.}, 228:282--293, 2009.

\bibitem{GLT07}
S.~G\H{u}nter, K.~Lackner, and C.~Tichmann.
\newblock Finite element and higher order difference formulations for modelling
  heat transport in magnetised plasmas.
\newblock {\em J. Comput. Phys.}, 226:2306--2316, 2007.

\bibitem{GYK05}
S.~G\H{u}nter, Q.~Yu, J.~Kruger, and K.~Lackner.
\newblock Modelling of heat transport in magnetised plasmas using non-aligned
  coordinates.
\newblock {\em J. Comput. Phys.}, 209:354--370, 2005.

\bibitem{NW00}
K.~Nishikawa and M.~Wakatani.
\newblock {\em Plasma Physics}.
\newblock Springer-Verlag Berlin Heidelberg, New York, 2000.

\bibitem{SH07}
P.~Sharma and G.~W. Hammett.
\newblock Preserving monotonicity in anisotropic diffusion.
\newblock {\em J. Comput. Phys.}, 227:123--142, 2007.

\bibitem{Sti92}
T.~Stix.
\newblock {\em Waves in Plasmas}.
\newblock Amer. Inst. Phys., New York, 1992.

\bibitem{ABBM98a}
I.~Aavatsmark, T.~Barkve, {\O}.~B{\o}e, and T.~Mannseth.
\newblock Discretization on unstructured grids for inhomogeneous, anisotropic
  media. {I}. {D}erivation of the methods.
\newblock {\em SIAM J. Sci. Comput.}, 19:1700--1716 (electronic), 1998.

\bibitem{ABBM98b}
I.~Aavatsmark, T.~Barkve, {\O}.~B{\o}e, and T.~Mannseth.
\newblock Discretization on unstructured grids for inhomogeneous, anisotropic
  media. {II}. {D}iscussion and numerical results.
\newblock {\em SIAM J. Sci. Comput.}, 19:1717--1736 (electronic), 1998.

\bibitem{CSW95}
P.~I. Crumpton, G.~J. Shaw, and A.~F. Ware.
\newblock Discretisation and multigrid solution of elliptic equations with
  mixed derivative terms and strongly discontinuous coefficients.
\newblock {\em J. Comput. Phys.}, 116:343--358, 1995.

\bibitem{EAK01}
T.~Ertekin, J.~H. Abou-Kassem, and G.~R. King.
\newblock {\em Basic Applied Reservoir Simulation}.
\newblock SPE textbook series, Vol. 7, Richardson, Texas, 2001.

\bibitem{MD06}
M.~J. Mlacnik and L.~J. Durlofsky.
\newblock Unstructured grid optimization for improved monotonicity of discrete
  solutions of elliptic equations with highly anisotropic coefficients.
\newblock {\em J. Comput. Phys.}, 216:337--361, 2006.

\bibitem{CS00}
T.~F. Chan and J.~Shen.
\newblock Non-texture inpainting by curvature driven diffusions ({CDD}).
\newblock {\em J. Vis. Commun. Image Rep}, 12:436--449, 2000.

\bibitem{CSV03}
T.~F. Chan, J.~Shen, and L.~Vese.
\newblock Variational {PDE} models in image processing.
\newblock {\em Not. AMS J.}, 50:14--26, 2003.

\bibitem{KM09}
D.~A. Karras and G.~B. Mertzios.
\newblock New {PDE}-based methods for image enhancement using {SOM} and
  {B}ayesian inference in various discretization schemes.
\newblock {\em Meas. Sci. Technol.}, 20:104012, 2009.

\bibitem{MS89}
D.~Mumford and J.~Shah.
\newblock Optimal approximations by piecewise smooth functions and associated
  variational problems.
\newblock {\em Commun. Pure Appl. Math}, 42:577--685, 1989.

\bibitem{PM90}
P.~Perona and J.~Malik.
\newblock Scale-space and edge detection using anisotropic diffusion.
\newblock {\em {IEEE} Trans. Pattern Anal. Mach. Intel.}, 12:629--639, 1990.

\bibitem{Wei98}
J.~Weickert.
\newblock {\em Anisotropic Diffusion in Image Processing}.
\newblock Teubner-Verlag, Stuttgart, Germany, 1998.

\bibitem{BKK08}
J.~Brandts, S.~Korotov, and M.~K\v{r}\'i\v{z}ek.
\newblock The discrete maximum principle for linear simplicial finite element
  approximations of a reaction-diffusion problem.
\newblock {\em Lin. Alg. Appl.}, 429:2344--2357, 2008.

\bibitem{Cia70}
P.~G. Ciarlet.
\newblock Discrete maximum principle for finite difference operators.
\newblock {\em Aequationes Math.}, 4:338--352, 1970.

\bibitem{CR73}
P.~G. Ciarlet and P.-A. Raviart.
\newblock Maximum principle and uniform convergence for thefinite element method.
\newblock {\em Comput. Meth. Appl. Mech. Engrg.}, 2:17--31, 1973.

\bibitem{DDS04}
A.~Dr\v{a}g\v{a}nescu, T.~F. Dupont, and L.~R. Scott.
\newblock Failure of the discrete maximum principle for an elliptic finite
  element problem.
\newblock {\em Math. Comp.}, 74:1--23, 2004.

\bibitem{Hua11}
W.~Huang.
\newblock Discrete maximum principle and a delaunay-type mesh condition for
  linear finite element approximations of two-dimensional anisotropic diffusion
  problems.
\newblock {\em Numer. Math. Theory Meth. Appl.}, 4:319--334, 2011.
\newblock (arXiv:1008.0562v1).

\bibitem{KK05}
J.~Kar{\'a}tson and S.~Korotov.
\newblock Discrete maximum principles for finite element solutions of nonlinear
  elliptic problems with mixed boundary conditions.
\newblock {\em Numer. Math.}, 99:669--698, 2005.

\bibitem{KKK07}
J.~Kar\'atson, S.~Korotov, and M.~K\v{r}\'i\v{z}ek.
\newblock On discrete maximum principles for nonlinear elliptic problems.
\newblock {\em Math. Comput. Sim.}, 76:99--108, 2007.

\bibitem{KSS09}
D.~Kuzmin, M.~J. Shashkov, and D.~Svyatskiy.
\newblock A constrained finite element method satisfying the discrete maximum
  principle for anisotropic diffusion problems.
\newblock {\em J. Comput. Phys.}, 228:3448--3463, 2009.

\bibitem{LH10}
X.~P. Li and W.~Huang.
\newblock An anisotropic mesh adaptation method for the finite element solution
  of heterogeneous anisotropic diffusion problems.
\newblock {\em J. Comput. Phys.}, 229:8072--8094, 2010 (arXiv:1003.4530v2).

\bibitem{LSS07}
X.~P. Li, D.~Svyatskiy, and M.~Shashkov.
\newblock Mesh adaptation and discrete maximum principle for {2D} anisotropic
  diffusion problems.
\newblock Technical Report LA-UR 10-01227, Los Alamos National Laboratory, Los
  Alamos, NM, 2007.

\bibitem{LS08}
R.~Liska and M.~Shashkov.
\newblock Enforcing the discrete maximum principle for linear finite element
  solutions of second-order elliptic problems.
\newblock {\em Comm. Comput. Phys.}, 3:852--877, 2008.

\bibitem{LuHuQi2012}
C.~Lu, W.~Huang, and J.~Qiu.
\newblock Maximum principle in linear finite element approximations of
  anisotropic diffusion-convection-reaction problems.
\newblock {\em (submitted)}, 2012 (arXiv:1201.3564v1).

\bibitem{ShYu2011}
Z.~Sheng and G.~Yuan.
\newblock The finite volume scheme preserving extremum principle for diffusion
  equations on polygonal meshes.
\newblock {\em J. Comput. Phys.}, 230:2588--2604, 2011.

\bibitem{Sto82}
G.~Stoyan.
\newblock On a maximum principle for matrices, and on conservation of
  monotonicity. {W}ith applications to discretization methods.
\newblock {\em Z. Angew. Math. Mech.}, 62:375--381, 1982.

\bibitem{Sto86}
G.~Stoyan.
\newblock On maximum principles for monotone matrices.
\newblock {\em Lin. Alg. Appl.}, 78:147--161, 1986.

\bibitem{SF73}
G.~Strang and G.~J. Fix.
\newblock {\em An Analysis of the Finite Element Method}.
\newblock Prentice Hall, Englewood Cliffs, NJ, 1973.

\bibitem{WaZh11}
J.~Wang and R.~Zhang.
\newblock Maximum principle for {P1}-conforming finite element approximations
  of quasi-linear second order elliptic equations.
\newblock 2011.
\newblock (arXiv:1105.1466v3).

\bibitem{YuSh2008}
C.~Yuan and Z.~Sheng.
\newblock Monotone finite volume schemes for diffusion equations on polygonal meshes.
\newblock {\em J. Comput. Phys.}, 227:6288--6312, 2008.

\bibitem{Das2006}
L.~Dascal.
\newblock Well-posedness and maximum principle for {PDE} based models in image processing.
\newblock PhD thesis, Tel-Aviv University, 2006.

{\Redink{\bibitem{Far2010}
I.~Farag{\'o}.
\newblock Discrete maximum principle for finite element parabolic models in higher dimensions.
\newblock {\em Math. Comput. Simulation}, 80:1601--1611, 2010.}}

\bibitem{FH06}
I.~Farag{\'o} and R.~Horv{\'a}th.
\newblock Discrete maximum principle and adequate discretizations of linear
  parabolic problems.
\newblock {\em SIAM J. Sci. Comput.}, 28:2313--2336, 2006.

\bibitem{FH07}
I.~Farag{\'o} and R.~Horv{\'a}th.
\newblock A review of reliable numerical models for three-dimensional linear
  parabolic problems.
\newblock {\em Int. J. Numer. Meth. Engng.}, 70:25--45, 2007.

\bibitem{FH09}
I.~Farag{\'o} and R.~Horv{\'a}th.
\newblock Continuous and discrete parabolic operators and their qualitative
  properties.
\newblock {\em IMA J. Numer. Anal.}, 29:606--631, 2009.

\bibitem{FHK05}
I.~Farag{\'o}, R.~Horv{\'a}th, and S.~Korotov.
\newblock Discrete maximum principle for linear parabolic problems solved on
  hybrid meshes.
\newblock {\em Appl. Numer. Math.}, 53:249--264, 2005.

\bibitem{FKK09}
I.~Farag{\'o}, Kar{\'a}tson, and S.~Korotov.
\newblock Discrete maximum principles for nonlinear parabolic pde systems.
\newblock Technical Report~93, Tampere University of Technology. Department of
  Mathematics, 2009.

\bibitem{Fuj1973}
H.~Fujii.
\newblock Some remarks on finite element analysis of time-dependent field
  problems.
\newblock In {\em Theory and Proactice in Finite Element Structural Analysis},
  pages 91--106. University of Tokyo, Tokyo, 1973.

\bibitem{Har04}
I.~Harari.
\newblock Stability of semidiscrete formulations for parabolic problems at
  small time steps.
\newblock {\em Comput. Methods Appl. Mech. Engrg.}, 193:1491--1516, 2004.

\bibitem{LE93}
M.~Lobo and A.~F. Emery.
\newblock The discrete maximum principle in finite-element thermal radiation
  analysis.
\newblock {\em Numer. Heat Transfer, Part B}, 24:209--227, 1993.

\bibitem{MVK89}
V.~Murti, S.~Valliappan, and N.~Khalili-Naghadeh.
\newblock Time step constraints in finite element analysis of the poisson type
  equations.
\newblock {\em Comput. Struct.}, 31:269--273, 1989.

\bibitem{TW08}
V.~Thom{\'e}e and L.~B. Wahlbin.
\newblock On the existence of maximum principles in parabolic finite element
  equations.
\newblock {\em Math. Comput.}, 77:11--19, 2008.

\bibitem{VKH08}
T.~Vejchodsk{\'y}, S.~Korotov, and A.~Hannukainen.
\newblock Discrete maximum principle for parabolic problems solved by prismatic
  finite elements.
\newblock Technical Report~77, Institute of Mathematics, AS CR, Prague, 2008.

\bibitem{YG06}
C.~Yang and Y.~Gu.
\newblock Minimum time-step criteria for the galerkin finite element methods
  applied to one-dimensional parabolic partial differential equations.
\newblock {\em Numer Meth. P. D. E.}, 22:259--273, 2006.

\bibitem{LePot05}
C.~Le~Potier.
\newblock Sch\'ema volumes finis monotone pour des op\'erateurs de diffusion
  fortement anisotropes sur des maillages de triangles non structur\'es.
\newblock {\em C. R. Math. Acad. Sci. Paris}, 341:787--792, 2005.

\bibitem{LePot09}
C.~Le~Potier.
\newblock A nonlinear finite volume scheme satisfying maximum and minimum principles for diffusion operators.
\newblock {\em Int. J. Finite Vol.}, 6(2):20 pp, 2009.

\bibitem{Cia78}
P.~G. Ciarlet.
\newblock {\em The Finite Element Method for Elliptic Problems}.
\newblock North-Holland, Amsterdam, 1978.

{\Redink{\bibitem{EN97}
J.~Emert and R.~Nelson.
\newblock Volume and Surface Area for Polyhedra and Polytopes.
\newblock {\em Math. Mag.} 70:365--371, 1997.}}

{\Redink{
\bibitem{HR11}
W.~Huang and R.~D.~Russell.
\newblock {\em Adaptive Moving Mesh Methods}.
\newblock Springer-Verlag Berlin Heidelberg, New York, 2011.
}}

\bibitem{Hec97}
F.~Hecht.
\newblock {BAMG -- Bidimensional Anisotropic Mesh Generator homepage}.
\newblock {http://www.ann.jussieu.fr/$\sim$hecht/ftp/bamg/}, 1997.

\bibitem{SCTPKB10}
A.~R. Sanderson, G.~Chen, X.~Tricoche, D.~Pugmire, S.~Kruger, and J.~Breslau.
\newblock Analysis of recurrent patterns in toroidal magnetic fields.
\newblock {\em IEEE Trans. Vis. Comput. Graph.}, 16:1431--1440, 2010.

\end{thebibliography}
\end{document}